\documentclass{amsproc}
\usepackage{breqn,float,amssymb,pdflscape,rotating,hyperref,xcolor,setspace,graphicx}
\makeatletter
\@namedef{subjclassname@2020}{%
  \textup{2020} Mathematics Subject Classification}
\makeatother
\newtheorem{theorem}{Theorem}[section]

\theoremstyle{definition}

\newtheorem{example}[theorem]{Example}

\theoremstyle{remark}

\numberwithin{equation}{section}
\allowdisplaybreaks
\begin{document}
\onehalfspacing
%
\title{A note on the Hurwitz-Lerch zeta function}


\author{Robert Reynolds}
\address[Robert Reynolds]{Department of Mathematics and Statistics, York University, Toronto, ON, Canada, M3J1P3}
\email[Corresponding author]{milver73@gmail.com}
\thanks{}




\date{}

\dedicatory{}

\begin{abstract}
In this work we derive a functional equation in terms of the Hurwitz-Lerch zeta function along with definite integrals in terms of the incomplete gamma and Hurwitz-Lerch zeta functions. The method used in these derivations is contour integration. Special cases in terms of fundamental constants are produced. 
\end{abstract}

\maketitle

\section{Introduction}
Functional equations involving the Hurwitz-Lerch zeta function, derived using definite integrals was first studied by Lipschitz \cite{lipschitz} in (1857). In his work he used certain definite integrals previously studied by Euler to analyse and derive series forms for the Hurwitz-Lerch zeta function. Malmsten \cite{malmsten} in (1867), used definite integrals involving nested logarithmic functions to derive functional equation involving the log-gamma and trigonometric functions. Lerch \cite{lerch} in (1877), studied an infinite series and subsequently derived a functional equation where the parameter range is over the complex plane. In 1899 Jonqui\'{e}re \cite{jonquiere} established a functional equation in terms of the polylogarithm and Hurwitz-zeta functions using contour integration and definite integrals. Fine \cite{fine} in (1951), established a generalized Riemann proof previously studied by Hurwitz to derive a functional equation for the Riemann zeta function. Apostol \cite{apostol} in (1951), derived a functional equation for the Hurwitz-Lerch zeta function based on the transformation theory of theta-functions. Special case evaluations of the Hurwitz-Lerch zeta function were also produced for negative integer values of the third parameter. Oberhettinger \cite{ober_l} in (1956),  showed that the application of Poisson's summation formula in its ordinary form to the more general case of Lerch's zeta function does not present the difficulties with respect to convergence which  arise in its special cases. Berndt \cite{berndt} in (1972),  produced two simple derivations of the functional equation $\sum\limits_{n=0}^{\infty}\frac{e^{2\pi inx}}{(n+a)^s}$ where the original proof is due to Lerch. The Wolfram functions website has a detailed list of Hurwitz-Lerch zeta functional equations see [Wolfram,\href{http://functions.wolfram.com/10.06.17.0009.01}{0.06.17.0009.01}]. Other resources for Hurwitz-Lerch zeta functional equations are Erd\'{e}yli et al. \cite{erd} pp. 27-31, Navas et al. \cite{navas}, Lagarias et al. \cite{lagarias} pp. 1-48, and [DLMF,\href{https://dlmf.nist.gov/25.14}{25.14}].
\\\\
In this paper we derive the functional equation given by
\begin{multline}\label{eq:fun_eq}
\Phi \left(e^{-2 i m \pi },-k,1-\frac{t}{2 \pi }\right)\\
=i (-1)^k e^{-\frac{1}{2} i (3 k \pi +2 m (-2 \pi +t))} (2 \pi )^{-1-k} \Gamma (1+k) \left(-\Phi \left(e^{-i t},1+k,m\right)\right. \\ \left.+(-1)^k e^{i t} \Phi \left(e^{i t},1+k,1-m\right)\right)
\end{multline}
where the parameters $k,t,m$ are general complex numbers and $Re(m)\geq 0$. This functional equation which has a wider range of evaluation for the parameters relative to previus results in the latter, will be used to simplify definite integrals which are expressed in terms of the Hurwitz-Lerch zeta function. The derivation involves two definite integrals and follows the method used by us in~\cite{reyn4}. This method involves using a form of the generalized Cauchy's integral formula given by
\begin{equation}\label{intro:cauchy}
\frac{y^k}{\Gamma(k+1)}=\frac{1}{2\pi i}\int_{C}\frac{e^{wy}}{w^{k+1}}dw.
\end{equation}
where $C$ is in general an open contour in the complex plane where the bilinear concomitant has the same value at the end points of the contour. We then multiply both sides by a function of $x$, then take a definite integral of both sides. This yields a definite integral in terms of a contour integral. Then we multiply both sides of Equation~(\ref{intro:cauchy})  by another function of $y$ and take the infinite sum of both sides such that the contour integral of both equations are the same.
\section{Contour integrals and definite integral representations}
in this section we will derive the contour integral representation for the definite integrals used to derive equation (\ref{eq:fun_eq}). 
\subsection{The first definite integral}
\begin{theorem}
For all $k,t,m\in\mathbb{C}$, where $Im(b)> 0,Re(m)<0$.
\begin{multline}\label{eq:first_int}
\int_0^{\infty } \frac{x^m \log ^k(a x)}{1-b x} \, dx\\
=-(-1)^m b^{-1-m} e^{i m \pi } (2 i \pi )^{1+k} \Phi \left(e^{2 i m
   \pi },-k,-\frac{i \left(i \pi +\log (a)+\log \left(-\frac{1}{b}\right)\right)}{2 \pi }\right)
\end{multline}
\end{theorem}
\begin{proof}
Use equation (7) in \cite{reyn5} with $\alpha=-b$ and $n=1$.
\end{proof}
\section{The second definite integral}
The contour integral representation for the addition of equations 6(9) and 6(10) in \cite{bdh},  is given by;
\begin{equation}\label{eq:sdi}
\frac{1}{2\pi i}\int_{0}^{1}\int_{C}\frac{a^w w^{-1-k} x^{-1+m+w}}{1-e^{i t} x}dwdx=\frac{1}{2\pi i}\int_{C}\sum_{n=0}^{\infty}\frac{a^w e^{i n t} w^{-1-k}}{m+n+w}dw
\end{equation}
 Using a generalization of Cauchy's integral formula \ref{intro:cauchy}, we form the definite integral by replacing $y$ by $\log{ax}$ and multiply both sides by $\frac{x^{m-1}}{1-e^{i t} x}$;
 \begin{multline}\label{eq:sdi1}
\int_{0}^{1}\frac{x^{-1+m} \log ^k(a x)}{\left(-1+e^{i t} x\right) k!}dx\\
=\frac{1}{2\pi i}\int_{0}^{1}\int_{C}\frac{w^{-1-k} x^{-1+m} (a x)^w}{-1+e^{i t} x}dwdx\\
=\frac{1}{2\pi i}\int_{C}\int_{0}^{1}\frac{w^{-1-k} x^{-1+m} (a x)^w}{-1+e^{i t} x}dwdx
\end{multline}
We are able to switch the order of integration over $x$ and $w$ using Fubini's theorem for multiple integrals see page 178 in \cite{gelca}, since the integrand is of bounded measure over the space $\mathbb{C} \times [0,1]$.
\subsubsection{The Incomplete Gamma~Function}
The incomplete gamma functions are given in equation [DLMF,\href{https://dlmf.nist.gov/8.4.E13}{8.4.13}], $\gamma(a,z)$ and $\Gamma(a,z)$, are defined by
\begin{equation}
\gamma(a,z)=\int_{0}^{z}t^{a-1}e^{-t}dt
\end{equation}
and
\begin{equation}
\Gamma(a,z)=\int_{z}^{\infty}t^{a-1}e^{-t}dt
\end{equation}
where $Re(a)>0$. The~incomplete gamma function has a recurrence relation given by
\begin{equation}
\gamma(a,z)+\Gamma(a,z)=\Gamma(a)
\end{equation}
where $a\neq 0,-1,-2,..$. The~incomplete gamma function is continued analytically by
\begin{equation}
\gamma(a,ze^{2m\pi i})=e^{2\pi mia}\gamma(a,z)
\end{equation}
and
\begin{equation}\label{eq:7}
\Gamma(a,ze^{2m\pi i})=e^{2\pi mia}\Gamma(a,z)+(1-e^{2\pi m i a})\Gamma(a)
\end{equation}
where $m\in\mathbb{Z}$, $\gamma^{*}(a,z)=\frac{z^{-a}}{\Gamma(a)}\gamma(a,z)$ is entire in $z$ and $a$. When $z\neq 0$, $\Gamma(a,z)$ is an entire function of $a$ and $\gamma(a,z)$ is meromorphic with simple poles at $a=-n$ for $n=0,1,2,...$ with residue $\frac{(-1)^n}{n!}$. These definitions are listed in [DLMF,\href{https://dlmf.nist.gov/8.2.i}{8.2(i)}] and [DLMF,\href{https://dlmf.nist.gov/8.2.ii}{8.2(ii)}].
The incomplete gamma functions are particular cases of the more general hypergeometric and Meijer G functions see section (5.6) and equation (6.9.2) in \cite{erd}. 
Some Meijer G representations we will use in this work are given by;
\begin{equation}\label{g1}
\Gamma (a,z)=\Gamma (a)-G_{1,2}^{1,1}\left(z\left|
\begin{array}{c}
 1 \\
 a,0 \\
\end{array}
\right.\right)
\end{equation}
and
\begin{equation}\label{g2}
\Gamma (a,z)=G_{1,2}^{2,0}\left(z\left|
\begin{array}{c}
 1 \\
 0,a \\
\end{array}
\right.\right)
\end{equation}
from equations (2.4) and (2.6a) in \cite{milgram}.  We will also use the derivative notation given by;
\begin{equation}\label{g3}
\frac{\partial \Gamma (a,z)}{\partial a}=\Gamma (a,z) \log (z)+G_{2,3}^{3,0}\left(z\left|
\begin{array}{c}
 1,1 \\
 0,0,a \\
\end{array}
\right.\right)
\end{equation}
from equations (2.19a) in \cite{milgram}, (9.31.3) in \cite{grad} and equations (5.11.1), (6.2.11.1) and (6.2.11.2) in \cite{luke}, and (6.36) in \cite{aslam}.
\subsubsection{Incomplete gamma function in terms of the contour integral}
In this section, we will once again use Cauchy's generalized integral formula, equation (\ref{intro:cauchy}), and take the infinite integral to derive equivalent sum representations for the contour integrals. We proceed using equation~(\ref{intro:cauchy}) and replace $y$ by $\log (a)+x$ and multiply both sides by $e^{m x}$ and simplify. Next, multiply both sides by $-e^{i n t}$ then take the infinite sum over $n\in[0,\infty)$ and simplify in terms of the incomplete gamma function to obtain
\begin{multline}\label{eq:sdi2}
-\sum_{n=0}^{\infty}\frac{a^{-m-n} e^{i n t} (-m-n)^{-1-k} \Gamma (1+k,-((m+n) \log (a)))}{k!}\\
=\frac{1}{2\pi i}\sum_{n=0}^{\infty}\int_{C}\frac{a^w e^{i n t} w^{-1-k}}{m+n+w}dw\\
=\frac{1}{2\pi i}\int_{C}\sum_{n=0}^{\infty}\frac{a^w e^{i n t} w^{-1-k}}{m+n+w}dw
\end{multline}
We are able to switch the order of integration and summation over $w$ using Tonellii's theorem for  integrals and sums see page 177 in \cite{gelca}, since the summand is of bounded measure over the space $\mathbb{C} \times [0,\infty)$.
%
\begin{theorem}
For all $k,a,m\in\mathbb{C},t\in\mathbb{R}$ then,
\begin{multline}\label{eq:second_int}
\int_0^1 \frac{x^{-1+m} \log ^k(a x)}{1-e^{i t} x} \, dx\\
=\sum _{n=0}^{\infty } a^{-m-n} e^{i n t} (-1)^k (m+n)^{-1-k}
   \Gamma (1+k,-((m+n) \log (a)))
\end{multline}
\end{theorem}
\begin{proof}
The right-hand sides of relations (\ref{eq:sdi1}) and (\ref{eq:sdi2}) are identical relative to equation (\ref{eq:sdi}), we can equate the left-hand sides. Simplify the gamma function yields the desired conclusion.
\end{proof}
\subsection{Derivation of the Hurwitz-Lerch zeta functional equation}
The derivation of equation (\ref{eq:fun_eq}) involves equations (\ref{eq:first_int}) and (\ref{eq:second_int}) where $a=1$. We first write equation (\ref{eq:first_int}) over $x\in[0,1]$ and $x\in[1,\infty)$. Next we transform the definite integral over $x\in[1,\infty)$ to $x\in[0,1]$ where we note the power of the variable $m$ is negative. Next we substitute these two integrals over $x\in[0,1]$ using equation (\ref{eq:second_int}) and simplify to yield the stated result.
\section{Alternate integral form over $x \in  [0,1] $}
In this section we will evaluate one of the definite integrals resulting from the subdivision in section (3.1). This integral is of interest due to its Cauchy principal evaluations. We will also use equation (\ref{eq:fun_eq}) to simplify the right-hand side whenever necessary.
\begin{theorem}
\begin{multline}
\int_0^1 \frac{x^{-1-m} \log ^k\left(\frac{a}{x}\right)}{1-b x} \, dx
=\sum _{n=0}^{\infty } \frac{(-1)^k \Gamma
   (1+k,-((1+m+n) \log (a)))}{b^{n+1} a^{1+m+n} (1+m+n)^{k+1}}\\
+(-1)^{2 m} b^m (2 i \pi )^{k+1} \Phi \left(e^{2 i m \pi
   },-k,-\frac{i (i \pi +\log (a)+\log (-b))}{2 \pi }\right)
\end{multline}
where $Im(b)< 0,Re(m)<0$.
\end{theorem}
\begin{example}
In this first example, we set $k=-1$ and replace $a$ by $e^{a+ib}$ and replace $b$ by $c$. Next we form a second equation by replacing $a$ by $-a$ and take the difference of the first and second equations and simplify.
\begin{multline}
\int_0^1 \frac{x^{-1-m}}{(-1+c x) \left(a^2+b^2+2 i b \log (x)-\log ^2(x)\right)} \, dx\\
=\frac{1}{2 a}\left(\sum _{n=0}^{\infty }
   c^{-1-n} e^{-((a+i b) (1+m+n))} \left(-e^{2 a (1+m+n)} E_1((a-i b) (1+m+n))\right.\right. \\ \left.\left.
+E_1(-((a+i b) (1+m+n)))\right)\right. \\ \left.
+(-1)^{2 m} c^m
   \left(-\Phi \left(e^{2 i m \pi },1,\frac{-i a+b+\pi -i \log (-c)}{2 \pi }\right)\right.\right. \\ \left.\left.
+\Phi \left(e^{2 i m \pi },1,\frac{i a+b+\pi -i
   \log (-c)}{2 \pi }\right)\right)\right)
\end{multline}
where $Im(c)< 0,Re(m)<0$.
\end{example}
\begin{example}
In this example we set $k=1$ and repeat the procedure in the previous example.
\begin{multline}
\int_0^1 \frac{x^{-1-m} \log \left(a^2+b^2+2 i b \log (x)-\log ^2(x)\right)}{1-c x} \, dx\\
=\sum _{n=0}^{\infty } \frac{c^{-1-n} }{1+m+n}\left(2 i \pi +e^{-((a+i b) (1+m+n))} \left(e^{2 a (1+m+n)} \Gamma
   (0,(a-i b) (1+m+n))\right.\right. \\ \left.\left.
+\Gamma (0,-((a+i b) (1+m+n)))\right)\right. \\ \left.
-2 \log (1+m+n)+\log ((a-i b) (1+m+n))+\log (-((a+i b) (1+m+n)))\right)\\
+\frac{4 c^m e^{2 i m \pi } \pi }{-1+e^{2 i m \pi }} \left(-i \log (2 i \pi )+e^{i m \pi
   } \sin (m \pi ) \left(\Phi'\left(e^{2 i m \pi },0,\frac{-i a+b+\pi -i \log (-c)}{2 \pi }\right)\right.\right. \\ \left.\left.
+\Phi'\left(e^{2 i m \pi },0,\frac{i a+b+\pi -i \log
   (-c)}{2 \pi }\right)\right)\right)-i c^m \pi  B_c(-m,0)
\end{multline}
where $Im(c)<0,Re(m)<0$.
\end{example}
\begin{example}
In this example we simply look at the case when $a=1$ and simplify the gamma function.
\begin{multline}\label{eq:44a}
\int_0^1 \frac{x^{-1-m} \log ^k\left(\frac{1}{x}\right)}{1-e^{-i t} x} \, dx\\
=-\exp (i t) \left(-(-1)^m e^{i m \pi } \left(e^{i t}\right)^{-1-m} (2 i \pi )^{1+k} \Phi \left(e^{2 i m \pi },-k,\frac{\pi
   -i \log \left(-e^{-i t}\right)}{2 \pi }\right)\right. \\ \left.-(-1)^k \Gamma (1+k) \Phi \left(e^{i t},1+k,1+m\right)\right)
\end{multline}
where $Re(m)<0$.
\end{example}
\begin{example}
In this example we look at the $n$-th derivative with respect to the parameter $b$ where $t=i\log(b)$.
\begin{multline}
\int_0^1 \frac{x^{-m}}{(b-x)^{n+1}} \, dx=\frac{(-1)^{-n} }{n!}\left(-b^{-m-n} \pi 
   (i+\cot (m \pi ))
   (-m)^{(n)}+\frac{\partial^n}{\partial b^n}\Phi(b,1,m)\right)
\end{multline}
where $Re(b)\leq 0$.
\end{example}
\begin{example}
In this example we use Mathematica by Wolfram to evaluate the previous definite integral and equate both expressions.
\begin{multline}
\frac{\partial^n}{\partial b^n}\Phi(b,1,m)\\
=b^{-m-n} \left(\frac{\pi  (i+\cot (m \pi )) \Gamma (1-m)}{\Gamma (1-m-n)}+(-1)^n B_{\frac{1}{b}}(1-m,-n) \Gamma (1+n)\right)
\end{multline}
where $0< Re(b)<1$.
\end{example}
\begin{example}
In this example we set $m=-1/2,t=\pi/2$ then take the first partial derivative with respect to $k$ and set $k=0$. Then we simplify the first partial derivative of the Hurwitz-Lerch zeta function using equation (\ref{eq:fun_eq}) and simplify.
\begin{multline}
\int_1^{\infty } \frac{\log (\log (x))}{\sqrt{x} (i+x)} \, dx\\
=\log \left(2^{-(-1)^{3/4} \left(\pi +2 i \log \left(\frac{15}{7} \cot \left(\frac{\pi }{8}\right)\right)\right)} \left(3-2
   \sqrt{2}\right)^{-\frac{1}{2} \sqrt[4]{-1} (\gamma +i \pi )} \right. \\ \left.
\exp \left(-\frac{\left(\frac{1}{2}-\frac{i}{2}\right) \left(\gamma  \left(\pi +4 i \tanh ^{-1}\left(\sqrt[4]{-1}\right)\right)+2 i \log
   \left(\frac{15}{7}\right) \log (4 \pi )\right)}{\sqrt{2}}\right)\right. \\ \left.
 \pi ^{-\frac{1}{2} (-1)^{3/4} \left(\pi +2 i \log \left(\frac{15}{7} \cot \left(\frac{\pi }{8}\right)\right)\right)} \left(\frac{3 \Gamma
   \left(-\frac{3}{8}\right)}{7 \Gamma \left(-\frac{7}{8}\right)}\right)^{\frac{1}{2} (-1)^{3/4} \pi } \left(\frac{5 \Gamma \left(-\frac{5}{8}\right)}{\Gamma \left(-\frac{1}{8}\right)}\right)^{\frac{3}{2}
   (-1)^{3/4} \pi }\right)\\
+\frac{1}{2} \sqrt[4]{-1}
   \left(-\zeta''\left(0,\frac{1}{8}\right)+\zeta''\left(0,\frac{3}{8}\right)+\zeta''\left(0,\frac{5}{8}\right)-\zeta''\left(0,\frac{7}{8}\right)\right)\
\end{multline}
\end{example}
\begin{example}
In this example we set $a=1,m=-1/2,t=\pi$ then take the first partial derivative with respect to $k$ then take the limit as $k\to 0$ using l'Hopital's rule and simplify in terms of the log-gamma function see [DLMF,\href{https://dlmf.nist.gov/25.11.E18}{25.11.18}] and Stieltjes constant [DLMF,\href{https://dlmf.nist.gov/25.2.E5}{25.2.5}].
\begin{multline}
\int_0^1 \frac{\log \left(\log \left(\frac{1}{x}\right)\right)}{\sqrt{x} (1+x)} \, dx\\
=\frac{1}{2} \left(\pi  \left(\gamma +\log \left(\frac{32 \pi ^2}{81}\right)-4 \text{log$\Gamma
   $}\left(-\frac{3}{4}\right)+4 \text{log$\Gamma $}\left(-\frac{1}{4}\right)\right)+\gamma _1\left(\frac{1}{4}\right)-\gamma _1\left(\frac{3}{4}\right)\right)
\end{multline}
\end{example}
\section{Special cases of the Hurwitz-Lerch zeta function giving simple constants}
In this section we look at special cases of equation (\ref{eq:fun_eq}) in terms of simple constants. Similar work can be found in the work by Guillera and Sondow \cite{guillera} and [Wolfram,\href{https://mathworld.wolfram.com/LerchTranscendent.html}{7-10}]
\begin{example}
In this example, set $m=1/2,t=\pi$, then take the first partial derivative with respect to $k$ then apply l'Hopital's rule as $k\to -1$ and simplify in terms of the log-gamma function see [DLMF,\href{https://dlmf.nist.gov/25.11.E18}{25.11.18}] and Stieltjes constant [DLMF,\href{https://dlmf.nist.gov/25.2.E5}{25.2.5}].
\begin{equation}
\gamma _1\left(\frac{1}{4}\right)-\gamma _1\left(\frac{3}{4}\right)=2 \pi  \log \left(\frac{3 e^{-\frac{\gamma }{2}} \Gamma \left(-\frac{3}{4}\right)}{2 \sqrt{2 \pi } \Gamma
   \left(-\frac{1}{4}\right)}\right)
\end{equation}
\end{example}
\begin{example}
In this example take the first partial derivative with respect to $k$ then set $m=1/2,k=1,t=0$ and simplify using equation [Wolfram,\href{https://mathworld.wolfram.com/LerchTranscendent.html}{8}].
\begin{equation}
\Phi'\left(1,2,\frac{1}{2}\right)=\frac{1}{2} \pi ^2 \log \left(\frac{4 \sqrt[3]{2} e^{\gamma } \pi }{A^{12}}\right)
\end{equation}
\end{example}
\begin{example}
In this example take the first partial derivative with respect to $k$ then set $m=1/2,t=0$ and take the limit as $k\to 2$ simplify using equation [Wolfram,\href{https://mathworld.wolfram.com/LerchTranscendent.html}{9}].
\begin{equation}
Li'_{-2}(-1)=-\frac{7 \zeta (3)}{4 \pi ^2}
\end{equation}
\end{example}
\begin{example}
In this example we simply set $t=0$ and simplify to get the polylogarithm is related to the Hurwitz zeta function by Jonqui\'{e}re's formula given by [Wolfram,\href{http://functions.wolfram.com/10.08.27.0004.01}{10.08.27.0004.01}].
\begin{multline}
\text{Li}_{-k}\left(e^{-2 i m \pi }\right)\\
=i (-1)^k e^{\frac{1}{2} (-3) i k \pi } (2 \pi )^{-1-k} \Gamma (1+k) \left((-1)^k \zeta (1+k,1-m)-\zeta (1+k,m)\right)
\end{multline}
\end{example}
\begin{example}
In this example we take the first partial derivative and set $m=1/2,k=0,t=\pi$ and simplify using equation [DLMF,\href{https://dlmf.nist.gov/25.11.vi}{25.11.18}].
\begin{equation}
\Phi'\left(-1,0,\frac{1}{2}\right)=\log \left(\frac{8 \Gamma \left(\frac{5}{4}\right)^2}{\pi }\right)
\end{equation}
\end{example}
\begin{example}
In this example we set $t=2\pi$, then take the first partial derivative with respect to $k$ then set $k=1,m=1/2$ and simplify using equation [Wolfram,\href{https://mathworld.wolfram.com/LerchTranscendent.html}{8}].
\begin{equation}
\Phi'(-1,-1,0)=\log \left(\frac{\sqrt[3]{2} \exp \left(\frac{1}{4}\right)}{A^3}\right)
\end{equation}
\end{example}
\begin{example}
In this example we take the first partial derivative with respect to $k$ then set $k=1,m=1/4$ and simplify using equation [Wolfram,\href{https://mathworld.wolfram.com/LerchTranscendent.html}{8}].
\begin{equation}
\Phi'(-i,-1,0)=\log \left(\frac{2^{7/6} \exp \left(\frac{1}{2}\right)}{\exp \left(\frac{2 i
   C}{\pi }\right) A^6}\right)
\end{equation}
\end{example}
\begin{example}
In this example we take the first partial derivative with respect to $k$ then set $k=1,m=3/4$ and simplify using equation [Wolfram,\href{https://mathworld.wolfram.com/LerchTranscendent.html}{8}].
\begin{equation}
\Phi'(i,-1,0)=\log \left(\frac{2^{7/6} \exp \left(\frac{2 i C}{\pi }\right) \exp
   \left(\frac{1}{2}\right)}{A^6}\right)
\end{equation}
\end{example}
\begin{example}
In this example we take the first partial derivative with respect to $k$ then set $k=2,m=1/2$ and simplify using equation [Wolfram,\href{https://mathworld.wolfram.com/LerchTranscendent.html}{9}]..
\begin{equation}
\Phi'(-1,-2,0)=-\frac{7 \zeta (3)}{4 \pi ^2}
\end{equation}
\end{example}
\begin{example}
In this example we set $t=-2\pi$, then take the limit using l'Hopital's rule as $k\to 1$ and $m=1/2$ and simplify using equation [Wolfram,\href{https://mathworld.wolfram.com/LerchTranscendent.html}{8}].
\begin{equation}
\Phi'(-1,-1,2)=\log \left(\frac{\exp \left(\frac{1}{4}\right) \sqrt[3]{2}}{A^3}\right)
\end{equation}
\end{example}
\begin{example}
In this example we set $t=4\pi$, then take the limit using l'Hopital's rule as $k\to 1$ and $m=1/2$ and simplify using equation [Wolfram,\href{https://mathworld.wolfram.com/LerchTranscendent.html}{8}].
\begin{equation}
\Phi'(-1,-1,-1)=-i \pi -\log \left(\frac{\sqrt[3]{2} \sqrt[4]{e}}{A^3}\right)
\end{equation}
\end{example}
\begin{example}
In this example we set $t=\pi,k->k-1$, then take the first partial derivative with respect to $k$ and take the limit using l'Hopital's rule as $k\to 0$ and simplify using equations [DLMF,\href{https://dlmf.nist.gov/25.11.vi}{25.11.18}], [Wolfram,\href{https://mathworld.wolfram.com/Euler-MascheroniConstant.html}{1}]
\begin{multline}
\frac{1}{2} \pi  \left(\gamma +\log \left(\frac{8 \pi ^3}{\Gamma \left(\frac{1}{4}\right)^4}\right)\right)
=\frac{1}{2} \left(3 \gamma  \pi +\pi  \left(2 i \pi +\log \left(\frac{400 \pi ^3}{3 \Gamma
   \left(\frac{1}{4}\right)}\right)+2 \text{log$\Gamma $}\left(-\frac{5}{4}\right)\right)\right. \\ \left.
+2 i \gamma  \left(\log \left(-\Gamma \left(-\frac{3}{4}\right)\right)-\text{log$\Gamma
   $}\left(-\frac{3}{4}\right)\right)-(\pi +2 i \log (4 \pi )) \text{log$\Gamma $}\left(-\frac{3}{4}\right)\right. \\ \left.
+i \left(12 \log ^2(2)-2 \log (3) \log (4 \pi )+2 \log (\pi ) \log \left(\Gamma
   \left(\frac{1}{4}\right)\right)+2 \log (4) \log \left(\pi  \Gamma
   \left(\frac{1}{4}\right)\right)\right.\right. \\ \left.\left.
-\zeta''\left(0,-\frac{1}{4}\right)+\zeta''\left(0,\frac{3}{4}\right)\right)\right)
\end{multline}
\end{example}
\begin{example}
In this example we set $t=2\pi(1-t),k=k-1,m=1/4$. Then take the first partial derivative with respect to $k$ and set $k=1$ and simplify using equation (6.11) in \cite{reyn_milli}.
\begin{multline}
-i \log \left(\frac{\Gamma \left(\frac{2+t}{4}\right)}{2 \Gamma \left(1+\frac{t}{4}\right)}\right)+\log \left(\frac{\Gamma \left(\frac{1+t}{4}\right)}{2 \Gamma
   \left(\frac{3+t}{4}\right)}\right)\\
=\frac{e^{\frac{1}{2} (-3) i \pi  (1+t)} }{6 \pi }\left(6 e^{2 i \pi  t} \, _2F_1\left(\frac{1}{4},1;\frac{5}{4};e^{2 i \pi  t}\right) (-2 i \gamma +\pi -2 i \log (2 \pi ))+2 \,
   _2F_1\left(\frac{3}{4},1;\frac{7}{4};e^{-2 i \pi  t}\right) \right. \\ \left.
\left(2 i \gamma +\pi +i \log \left(4 \pi ^2\right)\right)-3 i \Phi'\left(e^{-2 i \pi  t},1,\frac{3}{4}\right)+3 i
   e^{2 i \pi  t} \Phi'\left(e^{2 i \pi  t},1,\frac{1}{4}\right)\right)
\end{multline}
\end{example}
\begin{example}
In this example we set $t=2\pi(1-t),k=k-1,m=-1/4$. Then take the first partial derivative with respect to $k$ and set $k=1$ and simplify using equation (6.7) in \cite{reyn_aims}.
\begin{multline}
\log \left(\frac{\Gamma \left(\frac{t}{4}\right)}{2 \Gamma \left(\frac{2+t}{4}\right)}\right)+i \log \left(\frac{\Gamma \left(\frac{1+t}{4}\right)}{2 \Gamma
   \left(\frac{3+t}{4}\right)}\right)
=\frac{1}{6 \pi }\left(e^{\frac{1}{2} (-5) i \pi  t} \left(-3 i \Phi'\left(e^{-2 i \pi  t},1,\frac{5}{4}\right)\right.\right. \\ \left.\left.
+e^{2 i \pi  t} \left(-12 \pi +2 e^{2 i \pi  t}
   \, _2F_1\left(\frac{3}{4},1;\frac{7}{4};e^{2 i \pi  t}\right) (-2 i \gamma +\pi -2 i \log (2 \pi ))\right.\right.\right. \\ \left.\left.\left.
+6 \, _2F_1\left(\frac{1}{4},1;\frac{5}{4};e^{-2 i \pi  t}\right) (2 i \gamma +\pi +2 i \log (2 \pi ))+3
   i \Phi'\left(e^{2 i \pi  t},1,-\frac{1}{4}\right)\right)\right)\right)
\end{multline}
\end{example}
\begin{example}
In this example we look at functional identity where $Re(x) < 0$. A similar form is given in [Wolfram,\href{http://functions.wolfram.com/10.06.17.0009.01}{10.06.17.0009.01}] where $0<Re(x)<1$.
\begin{multline}\label{eq:515}
\Phi \left(e^{2 i \pi  x},1-s,a\right)\\
=-e^{-\frac{1}{2} i \pi  (-3+s+4 a (1+x))} (2 \pi )^{-s} \Gamma (s) \left(e^{i \pi  s} \Phi \left(e^{-2 i a \pi },s,1+x\right)+e^{2 i a \pi } \Phi \left(e^{2 i a
   \pi },s,-x\right)\right)
\end{multline}
\end{example}
\begin{example}
In this example using equation (\ref{eq:515}), we set $x=-1/2$ then take the first partial derivative with respect to $s$ then set $s=1$ and simplify using [DLMF,\href{https://dlmf.nist.gov/25.11.vi}{25.11.18}].
\begin{multline}
\frac{\log (2)}{2}-\log (-2+a)+\log (-1+a)-\text{log$\Gamma $}\left(-1+\frac{a}{2}\right)+\text{log$\Gamma $}\left(\frac{1}{2} (-1+a)\right)\\
=-\frac{1}{2 \pi }\left(e^{-i a \pi }
   \left(\Phi'\left(e^{-2 i a \pi },1,\frac{1}{2}\right)+e^{i a \pi } \left(\tanh ^{-1}\left(e^{i a \pi }\right) (2 \gamma +i \pi +\log (4)+2 \log (\pi ))\right.\right.\right. \\ \left.\left.\left.
+\coth ^{-1}\left(e^{i a
   \pi }\right) (-2 \gamma +i \pi -2 \log (2 \pi ))-e^{i a \pi } \Phi'\left(e^{2 i a \pi },1,\frac{1}{2}\right)\right)\right)\right)
\end{multline}
\end{example}
\section{A few evaluations}
\begin{example}
Using equation (\ref{eq:44a}) where $t\to i\log(b),m\to m-1$, we take the first partial derivative with respect to $k$ and set $m=1/2,k=1,b=-1$ in terms of Catalan's constant ($C$) and simplify using equation (\ref{eq:fun_eq}) and [Wolfram,\href{https://mathworld.wolfram.com/LerchTranscendent.html}{10}].
\begin{multline}
\int_0^1 \frac{\log \left(\frac{1}{x}\right) \log \left(\log \left(\frac{1}{x}\right)\right)}{\sqrt{x} (1+x)} \, dx=4 \left(C \log (4 \pi )+\pi 
   \left(-\zeta''\left(-1,\frac{1}{4}\right)+\zeta'\left(-1,\frac{3}{4}\right)\right)\right)
\end{multline}
\end{example}
\begin{example}
Using equation (\ref{eq:44a}) where $t\to i\log(b),m\to m-1$, we take the first partial derivative with respect to $k$ and set $k=1,b=1$ and $m=-1/2$ for the first equation and $m=1/2$ for the second equation and take their difference.
\begin{multline}
\int_0^1 \frac{\log \left(\frac{1}{x}\right) \log \left(\log \left(\frac{1}{x}\right)\right)}{\sqrt{x}} \, dx
=4-4 \gamma +4 i \pi
   +\Phi'\left(1,2,-\frac{1}{2}\right)-\Phi'\left(1,2,\frac{1}{2}\right)
\end{multline}
\end{example}
\begin{example}
Using equation (\ref{eq:44a}) where $t\to i\log(b),m\to m-1$, we take the first partial derivative with respect to $k$ and set $k=2,b=1$ and $m=-1/2$ for the first equation and $m=1/2$ for the second equation and take their difference.
\begin{multline}
\int_0^1 \frac{\log ^2\left(\frac{1}{x}\right) \log \left(\log \left(\frac{1}{x}\right)\right)}{\sqrt{x}} \, dx
=2 \left(12-8 \gamma +8 i \pi
   -\Phi'\left(1,3,-\frac{1}{2}\right)+\Phi'\left(1,3,\frac{1}{2}\right)\right)
\end{multline}
\end{example}
\begin{example}
In this example we using equation (\ref{eq:44a}) we form a second equation by replacing $m\to s$ then take their difference. Next we set $m=1/2,s=-1/2,b=2$ then take the first partial derivative with respect to $k$ and use l'Hopitals' rule as $k\to 0$ and simplify using [Wolfram,\href{https://mathworld.wolfram.com/Euler-MascheroniConstant.html}{1}] and [DLMF,\href{https://dlmf.nist.gov/25.11.E18}{25.11.18}]. Note the Cauchy principal value of the integral, where the is an essential singularity is at $x=1/2$. The range of integration is $x\in[0,-i]\cup [-i,1]$.
\begin{multline}
\int_0^1 \frac{(x-1) \log \left(\log \left(\frac{1}{x}\right)\right)}{\sqrt{x} (-1+2 x)} \, dx\\
=\frac{1}{8} \left(3 \sqrt{2} \pi ^2-2 \pi  \left(4 i+\sqrt{2} \left(\pi +i \left(\log (\pi )\right.\right.\right.\right. \\ \left.\left.\left.\left.
-2
   \left(\log \left(-2 i \sqrt{2} \pi +\log (2)+\sqrt{2 \left(-2 \pi ^2-2 i \sqrt{2} \pi  \log (2)+\log ^2(2)\right)}\right)\right.\right.\right.\right.\right. \\ \left.\left.\left.\left.\left.
-\log (\log (4))-\text{log$\Gamma $}\left(-\frac{i \log (2)}{4 \pi
   }\right)+\text{log$\Gamma $}\left(-\frac{1}{2}-\frac{i \log (2)}{4 \pi }\right)\right)\right)\right)\right)\right. \\ \left.
+4 \left(-\gamma  \left(2+\sqrt{2} \sinh ^{-1}(1)\right)+\log
   (16)+\Phi'\left(\frac{1}{2},1,-\frac{1}{2}\right)\right)\right)
\end{multline}
\end{example}
\begin{figure}[h]
\caption{Plot of $f(x)=Re\left(\frac{(x-1) \log \left(\log \left(\frac{1}{x}\right)\right)}{\sqrt{x} (2 x-1)}\right)$}
\centering
\includegraphics[width=0.9\textwidth]{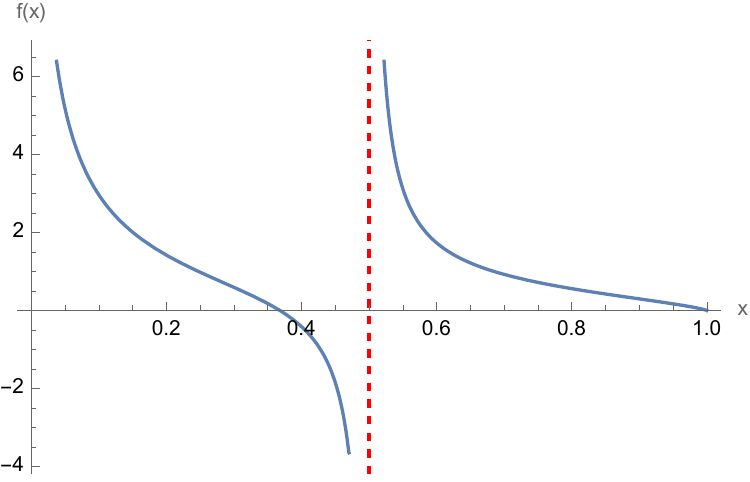}
\end{figure}
\begin{example}
In this example we form two equations using (\ref{eq:second_int}) with $k=-1$ and $a\to e^{ai}$ and  $a\to e^{-ai}$, take their difference and simplify.
\begin{multline}\label{eq:6.5}
\int_0^1 \frac{x^{m-1}}{(b x-1) \left(a^2+\log ^2(x)\right)} \,
   dx\\
=\frac{i  }{2 a}\sum _{n=0}^{\infty } b^n\left(e^{-i a (m+n)} \Gamma (0,-i a
   (m+n))-e^{i a (m+n)} \Gamma (0,i a (m+n))\right)
\end{multline}
where $Re(b)\leq 0$.
\end{example}
\begin{example}
In this example we use equation  (\ref{eq:second_int}) and take the $k$-th partial derivative with respect to $b$ and simplify.
\begin{multline}
\int_0^1 \frac{x^{-1+k+m}}{\left(a^2+\log ^2(x)\right) (-1+b x)^{k+1}}
   \, dx\\
=\sum _{n=0}^{\infty } \frac{i (-1)^{-k} b^{-k+n} \Gamma (1+n)
   \left(e^{-i a (m+n)} \Gamma (0,-i a (m+n))-e^{i a (m+n)} \Gamma (0,i a
   (m+n))\right)}{2 a \Gamma (1+k) \Gamma (1-k+n)}
\end{multline}
where $Re(b)\ <1$.
\end{example}
\begin{example}
In this example we use equation (\ref{eq:second_int}) and replace $t\to -i\log(b)$. Next we take the $q$-th derivative with respect to $b$ and replace $k\to s,q\to k$ and simplify.
\begin{multline}
\int_0^1 \frac{x^{-1+k+m} \log ^s(a x)}{(1-b x)^{k+1}} \, dx\\
=\sum
   _{n=0}^{\infty } \frac{(-1)^s a^{-m-n} b^{-k+n} (m+n)^{-1-s} \Gamma (1+n)
   \Gamma (1+s,-((m+n) \log (a)))}{\Gamma (1+k) \Gamma (1-k+n)}
\end{multline}
where $Re(m)> 0$.
\end{example}
\begin{example}
Generalized forms for Table 327 in \cite{bdh}. This example is derived using equation (\ref{eq:second_int}). First replace $t\to -i \log{b}$ then replace $b\to b+i c$. Next form a second equation by replacing $c\to -c$ and take their difference. We keep $|Re(b)|<1$ in order for the series to converge.
\begin{multline}\label{eq:6.8}
\int_0^1 \frac{x^m \log ^k(a x)}{1-2 b x+\left(b^2+c^2\right) x^2} \, dx\\
=\sum _{n=0}^{\infty } \frac{i (-1)^k
   \left((b-i c)^n-(b+i c)^n\right) \Gamma (1+k,-((m+n) \log (a)))}{2 c a^{m+n} (m+n)^{k+1}}
\end{multline}
where $0< Re(b) <1,0< Re(c) <1,0< Re(m) <1$.
\end{example}
\begin{example}
This example is derived using equation (\ref{eq:6.8}) by setting $k=-1,a=e^{ai}$. Next form a second equation by replacing $a\to -a$ and take the difference.
\begin{multline}\label{eq:}
\int_0^1 \frac{x^m}{\left(1-2 b x+\left(b^2+c^2\right) x^2\right) \left(a^2+\log ^2(x)\right)} \, dx\\
=\sum
   _{n=0}^{\infty } \frac{\left((b-i c)^n-(b+i c)^n\right) \left(\Gamma (0,-i a (m+n))-e^{2 i a (m+n)} \Gamma (0,i a
   (m+n))\right)}{4 a c e^{i a (m+n)}}
\end{multline}
where $0< Re(b) <1,0< Re(c) <1,0< Re(m) <1$.
\end{example}
\begin{example}
This example is derived using equation (\ref{eq:6.8}) and applying l'Hopital's rule as $c\to 0$ and simplify.
\begin{equation}\label{eq:6.10}
\int_0^1 \frac{x^m \log ^k(a x)}{(-1+b x)^2} \, dx=\sum _{n=0}^{\infty } \frac{(-1)^k b^{-1+n} n \Gamma (1+k,-((m+n)
   \log (a)))}{a^{m+n} (m+n)^{k+1}}
\end{equation}
where $0< Re(b) <1,0< Re(m) <1$.
\end{example}
\begin{example}
This example is derived by using equation (\ref{eq:6.10}) and taking the $k$-th partial derivative. We simplify using the Pochammer symbol in [DLMF,\href{https://dlmf.nist.gov/5.2.iii}{5.2.5}]. We also replace $n\to j,k\to n$.
\begin{equation}\label{eq:6.11}
\int_0^1 \frac{x^n \left(x^s-x^m\right)}{(b x-1)^{n+2} \log (x)} \, dx=\sum _{j=0}^{\infty } \frac{b^{-1+j-n}
   (j-n)_{n+1} \log \left(\frac{j+s}{j+m}\right)}{(-1)^n \Gamma (n+2)}
\end{equation}
where $0< Re(b) <1,0< Re(m) <1,0< Re(s) <1,n\in\mathbb{Z_{+}}$.
\end{example}
\begin{example}
A special case of equation (\ref{eq:6.11}) with $b=1/2,s=1/2,m=-1/2,n=0$.
\begin{equation}
\int_0^1 \frac{x-1}{(x-2)^2 \sqrt{x} \log (x)} \, dx=\sum _{j=0}^{\infty } 2^{-1-j} j \log \left(\frac{1+2
   j}{-1+2 j}\right)
\end{equation}
\end{example}
\begin{example}
In this example we look at the definite integral involving the nested logarithmic function originally studied by Malmsten \cite{malmsten}. We use equation (\ref{eq:6.10}), and set $a=1$, then apply l'Hopital's rule as $c\to 0$. Next we take the $j$-th partial derivative with respect to $b$ and simplify using the Pochammer symbol in [DLMF,\href{https://dlmf.nist.gov/5.2.iii}{5.2.5}]. We then take the first partial derivative with respect $k$ and set $k=0$. We then make a change of variables $n\to j, k\to s$ and take the sum over the range $j\in [0,\infty)$ since the series is zero at $j=0$.
\begin{multline}
\int_0^1 \frac{x^{s-1} \log (\log (x))}{(b x-1)^n} \, dx\\
=\frac{(-b)^{1-n}
  }{\Gamma (n)} \sum _{j=1}^{\infty } \frac{b^j j (\gamma -i \pi +\log (1+j-n+s))
   (2+j-n)_{n-2}}{1+j-n+s}
\end{multline}
where $|Re(b)|<1,0< Re(m)<1,0< Re(s)<1,n\in\mathbb{Z_{+}}$.
\end{example}
\begin{example}
In this example we look at the $n$-th partial derivative of the Hurwitz-Lerch zeta function expressed as a series. We use equations (\ref{eq:second_int}) and (\ref{eq:6.10}) with $a=1$. Using equation (\ref{eq:6.10}), we simply take the $j$-th partial derivative with respect to $b$ and algebraically made the left-hand sides of both equations the same.
\begin{equation}
\sum _{j=0}^{\infty } \frac{j z^{j-n}
   (1+j-n)_{-1+n}}{(a+j)^s}=\Phi^n(z,s,a)
\end{equation}
where $|Re(z)|<1,0< Re(a)<1,0< Re(s)<1,n\in\mathbb{Z_{+}}$.
\end{example}
\begin{example}
In this example we apply simple algebraic techniques to equation (\ref{eq:second_int}). We start by replacing $b\to -b$ to form a second equation, then take their difference. Next set $k=-1$ and replace $a\to -a$ to form another equation and take their difference and replace $a\to ai,b\to bi$ and simplify. Similar forms are listed in sections (4.282-4) in \cite{grad} and sections (2.6.18) in \cite{prud1}. This is an extended form of equations (2.6.18.10-11) in \cite{prud1}.
\begin{multline}
\int_0^1 \frac{x^{m-1}}{\left(b^2+x^2\right) \left(a^2+\log ^2\left(\frac{1}{x}\right)\right)} \, dx\\
=\sum
   _{n=0}^{\infty } \frac{i^{2 n+1} \cos \left(\frac{n \pi }{2}\right) \left(e^{i a (m+n)} \Gamma (0,i a
   (m+n))-e^{-i a (m+n)} \Gamma (0,-i a (m+n))\right)}{2 a b^{n+2}}
\end{multline}
where $0< Re(m)<1$.
\end{example}
\begin{example}
Extended version of equation (4.327.1) and  (4.327.2) in \cite{grad}. We use  equation (\ref{eq:second_int}) and take the first partial derivative with respect to $k$ then set $k=0$. Next we replace $a\to -a$ to form a second equation then add the two equations and simplify.
\begin{multline}
\int_0^1 \frac{x^{-1+m} \log \left(a^2+\log ^2\left(\frac{1}{x}\right)\right)}{b^2+x^2} \, dx\\
=\sum
   _{n=0}^{\infty } \frac{i^n \left(1+(-1)^n\right) \left(\Gamma (0,-i a (m+n))+e^{2 i a (m+n)} \Gamma (0,i a (m+n))+2
   e^{i a (m+n)} \log (a)\right)}{2 (m+n) b^{n+2} e^{i a (m+n)}}
\end{multline}
where $0< Re(m)<1,|Re(b)|>1$.
\end{example}
\begin{example}
Derivation of entries (2.6.20.1-4) in \cite{prud1}. We use  equation (\ref{eq:second_int}) set $a=0,m=1$ and replace $b\to -e^{i\lambda}$, then take the first partial derivative respect to $k$ and set $k=0$ and simplify.
\begin{equation}
\int_0^1 \frac{\log \left(\log \left(\frac{1}{x}\right)\right)}{e^{i \lambda }+x} \, dx=-\gamma  \log
   \left(1+e^{-i \lambda }\right)Li'_{1}\left(-e^{-i \lambda }\right)
\end{equation}
where $\lambda\in\mathbb{C}$.
\end{example}
\begin{example}
Derivation of generalized equation (2.6.20.7) in \cite{prud1} in terms of the Hurwitz-Lerch zeta function. We use  equation (\ref{eq:second_int}) and set $a=0$. Then form a second equation by replacing $b\to -b$ and take their difference. Next replace $b\to bi$, then take the first partial derivative with respect to $k$ and set $k=-1/2$ and simplify. 
\begin{multline}
\int_0^1 \frac{x^{-1+m} \log \left(\log \left(\frac{1}{x}\right)\right)}{\left(b^2+x^2\right) \sqrt{\log
   \left(\frac{1}{x}\right)}} \, dx\\
=\frac{\sqrt{\pi } \left(-\left(\left(\Phi
   \left(-\frac{i}{b},\frac{1}{2},m\right)+\Phi \left(\frac{i}{b},\frac{1}{2},m\right)\right) (\gamma +\log
   (4))\right)+\Phi'\left(-\frac{i}{b},\frac{1}{2},m\right)+\Phi'\left(\frac{i}{b},\frac{1}{2},m\right)\right)}{2 b^2}
\end{multline}
where $b,m\in\mathbb{C}$.
\end{example}
\begin{example}
Derivation of equation (2.6.20.4) in \cite{prud1}.  We use  equation (\ref{eq:second_int}) and set $a=0,m=1$ and simplify in terms of the Polylogarithm function using [DLMF,\href{https://dlmf.nist.gov/25.14.E2}{25.14.3}], then form a second equation replacing $b\to c$ and take their difference. Next replace $b\to e^{i\gamma},c\to e^{-i\gamma}$. Next we simplify the Polylogarithm function using [DLMF,\href{https://dlmf.nist.gov/25.12.E13}{25.12.13 }]. Next we take the first partial derivative with respect to $k$ and apply l'Hopitals rule as $k\to 0$ and simplify. intermediate steps are included.
\begin{multline}
\int_0^1 \frac{\log \left(\log \left(\frac{1}{x}\right)\right)}{1+x^2-2 x \cos (\gamma )} \, dx
=-\frac{1}{4} i
   \csc (\gamma ) \left(\log (64) \log (\pi )+\left(\log \left(-e^{-i \gamma }\right)-\log \left(-64 e^{i \gamma
   }\right)\right) \log (2 \pi )\right. \\ \left.
+i \left(\pi +i \log \left(4 \pi ^2\right)\right) \log \left(-\pi -i \log \left(-e^{-i
   \gamma }\right)\right)+\left(-i \pi +\log \left(\frac{1}{4 \pi ^2}\right)\right) \log \left(-\pi +i \log
   \left(-e^{-i \gamma }\right)\right)\right. \\ \left.
+\left(-i \pi +\log \left(4 \pi ^2\right)\right) \log \left(-\pi -i \log
   \left(-e^{i \gamma }\right)\right)+(i \pi +2 \log (\pi )) \log \left(-\pi +i \log \left(-e^{i \gamma
   }\right)\right)\right. \\ \left.
+\log (4) \log \left(-8 \pi +8 i \log \left(-e^{i \gamma }\right)\right)+i \left((1+2 i) \pi +\log
   \left(4^i \pi ^{2 i}\right)\right) \text{log$\Gamma $}\left(-\frac{\pi +i \left(\log (-1)+\log \left(e^{-i \gamma
   }\right)\right)}{2 \pi }\right)\right. \\ \left.
-\left((2+i) \pi +i \log \left(\left(\frac{1}{4}\right)^i \pi ^{-2 i}\right)\right)
   \text{log$\Gamma $}\left(-\frac{1}{2}+\frac{i \left(\log (-1)+\log \left(e^{-i \gamma }\right)\right)}{2 \pi
   }\right)\right. \\ \left.
-i \left((1+2 i) \pi +\log \left(4^i \pi ^{2 i}\right)\right) \text{log$\Gamma $}\left(-\frac{\pi +i
   \left(-\log (-1)+\log \left(e^{i \gamma }\right)\right)}{2 \pi }\right)\right. \\ \left.
+\left((2+i) \pi +i \log \left(4^{-i} \pi
   ^{-2 i}\right)\right) \text{log$\Gamma $}\left(-\frac{1}{2}+\frac{i \left(-\log (-1)+\log \left(e^{i \gamma
   }\right)\right)}{2 \pi }\right)\right. \\ \left.
+\zeta''\left(0,\frac{\pi -i \left(\log (-1)+\log \left(e^{-i
   \gamma }\right)\right)}{2 \pi }\right)
+\zeta''\left(0,\frac{\pi +i \left(\log (-1)+\log
   \left(e^{-i \gamma }\right)\right)}{2 \pi }\right)\right. \\ \left.
-\zeta''\left(0,\frac{\pi -i \left(-\log
   (-1)+\log \left(e^{i \gamma }\right)\right)}{2 \pi }\right)
-\zeta''\left(0,\frac{\pi +i
   \left(-\log (-1)+\log \left(e^{i \gamma }\right)\right)}{2 \pi }\right)\right)\\
=-\frac{1}{2} \csc (\gamma )
   \left(\gamma  \log (2 \pi )+\pi  \log \left(-\frac{1}{2 \pi }+\frac{1}{\gamma }\right)-\pi  \text{log$\Gamma
   $}\left(-\frac{\gamma }{2 \pi }\right)+\pi  \text{log$\Gamma $}\left(-1+\frac{\gamma }{2 \pi
   }\right)\right)\\
=-\frac{1}{2} \pi  \csc (\gamma ) \log \left(\frac{(2 \pi
   )^{-1+\frac{\gamma }{\pi }} \Gamma \left(\frac{\gamma }{2 \pi }\right)}{\Gamma \left(1-\frac{\gamma }{2 \pi
   }\right)}\right)
\end{multline}
where $\gamma\in\mathbb{C}$.
\end{example}
\begin{example}
Derivation of a generalized form of equation (2.6.20.5) in \cite{prud1}. We use  equation (\ref{eq:second_int}) and set $a=0$ and simplify in terms of the Hurwitz-Lerch zeta function. Next we form a second equation by replacing $b\to -b$ and take their difference and simplify. Next we replace $x\to u^{n+1}$ and simplify. We then form another equation by replacing $m\to s$ and take the difference and simplify. Next we set $b=1$ and simplify the right-hand side in terms of the Hurwitz zeta function using equation [DLMF,\href{https://dlmf.nist.gov/25.14.E2}{25.14.2}] and equation (64:13:3) in \cite{atlas}. Next we replace $m \to n/(1 + n), s \to (2 + n)/(1 + n)$ and take the first partial derivative with respect to $k$. We then applyl 'Hopital's rule as $k\to 0$ and simplify interms of the Euler's constant $\gamma$ from [DLMF,\href{https://dlmf.nist.gov/5.2.E3}{5.2.3}] and Stieltjes constant $\gamma_{n}$ from [DLMF,\href{https://dlmf.nist.gov/25.2.E5}{25.2.5}] and simplify.
\begin{multline}
\int_0^1 \frac{x^{-1+z} \left(1-x^2\right) \log \left(\log \left(\frac{1}{x}\right)\right)}{-1+x^{2+2 z}} \,
   dx\\
=\frac{1}{4 (1+z)}\left(-((1+z) \log (4))+2 \pi  \left(\cot \left(\frac{\pi  z}{1+z}\right) (\gamma +\log (1+z))\right.\right. \\ \left.\left.
+\csc
   \left(\frac{\pi }{1+z}\right) (\gamma +\log (2+2 z))\right)+2 \gamma _1\left(\frac{z}{1+z}\right)+\gamma
   _1\left(\frac{z}{2+2 z}\right)-\gamma _1\left(\frac{2+z}{2+2 z}\right)\right. \\ \left.
-\gamma _1\left(1-\frac{1}{2 (1+z)}\right)-2
   \gamma _1\left(1+\frac{1}{1+z}\right)+\gamma _1\left(1+\frac{1}{2+2 z}\right)\right)
\end{multline}
where $z\in\mathbb{C}$.
\end{example}
\begin{example}
In this example, we use  equation (\ref{eq:second_int}) and take the $q$-th derivative of $m$ and simplify.
\begin{multline}\label{eq:6.21}
\int_0^1 \frac{x^{-1+m} \log ^k\left(\frac{1}{x}\right) \log ^q\left(\frac{1}{x}\right)}{-b+x} \,
   dx=-\frac{\Gamma (1+k) \Phi \left(\frac{1}{b},1+k+q,m\right) (1+k)_q}{b}
\end{multline}
where $|Re(b)|>1$.
\end{example}
\begin{example}
In this example we use equation (\ref{eq:6.21}) and take the first partial derivative with respect to $k$ and set $k=0$. Next form a second equation by replacing $b\to -b$ and take their difference and simplify.
\begin{multline}
\int_0^1 \frac{x^{-1+m} \log ^q\left(\frac{1}{x}\right) \log \left(\log
   \left(\frac{1}{x}\right)\right)}{b^2-x^2} \, dx\\
=\frac{\Gamma (1+q) }{2 b^2}\left(\left(\Phi
   \left(-\frac{1}{b},1+q,m\right)+\Phi \left(\frac{1}{b},1+q,m\right)\right) \psi
   ^{(0)}(1+q)\right. \\ \left.
+\Phi'\left(-\frac{1}{b},1+q,m\right)+\Phi'\left(\frac{1}{b},1+q,m\right)\right)
\end{multline}
where $|Re(b)|>1$.
\end{example}
\begin{example}
In this example, we use  equation (\ref{eq:second_int}) set $a=1$ and take the $p$-th derivative of $m$ and the first partial derivative with respect to $k$ and set $k=0$ simplify.
\begin{multline}
\int_0^1 \frac{x^{m+p-1} \log \left(\log \left(\frac{1}{x}\right)\right)}{(b x-1)^{p+1}} \, dx=\frac{(-1)^{-p}
  }{\Gamma(1+p)} \left(\gamma  \frac{\partial}{\partial b^p}\Phi^p(b,1,m)-\frac{\partial}{\partial b^p}\Phi^p(b,1,m)\right)
\end{multline}
where $|Re(b)|>1$.
\end{example}
\begin{example}
In this example, we use  equation (\ref{eq:second_int}) set $a=1,b=-1,k=2,m=1/2$ and simplify using equation [Wolfram,\href{https://mathworld.wolfram.com/LerchTranscendent.html}{7}].
\begin{equation}
\int_0^1 \frac{\log \left(\frac{1}{x}\right)}{(1+x) \sqrt{x}} \, dx=4 C
\end{equation}
\end{example}
\begin{example}
In this example, we use  equation (\ref{eq:second_int}) take the first partial derivative with respect to $k$ set $a=1,b=-1,k=1,m=1/2$ and simplify using equation [Wolfram,\href{https://mathworld.wolfram.com/LerchTranscendent.html}{8}].
\begin{equation}
\int_0^1 \frac{\log \left(\log \left(\frac{1}{x}\right)\right)}{\sqrt{x} (1+x)} \, dx=\frac{1}{2} \pi  \log
   \left(\frac{8 \pi ^3}{\Gamma \left(\frac{1}{4}\right)^4}\right)
\end{equation}
\end{example}
\begin{example}
In this example, we use  equation (\ref{eq:second_int}) take the first partial derivative with respect to $k$ set $a=1,b=-1,k=1,m=1$ and simplify using equation [Wolfram,\href{https://mathworld.wolfram.com/LerchTranscendent.html}{8}].
\begin{equation}
\int_0^1 \frac{\log \left(\log \left(\frac{1}{x}\right)\right)}{1+x} \, dx=-\frac{1}{2} \log ^2(2)
\end{equation}
\end{example}
\begin{example}
In this example we derive an integral representation for the inverse tangent integral. We use  equation (\ref{eq:second_int}) and set $a=0,b = -z^2, k = s, m = 1/2$ and simplify using equation [Wolfram,\href{https://mathworld.wolfram.com/InverseTangentIntegral.html}{6}].
\begin{equation}
\int_0^1 \frac{\log ^{-1+s}\left(\frac{1}{x}\right)}{\sqrt{x} \left(1+x z^2\right)} \, dx=\Gamma (s) \Phi
   \left(-z^2,s,\frac{1}{2}\right)
\end{equation}
where $s,z\in\mathbb{C}$.
\end{example}
\begin{example}
In this example we derive an integral representation for the Legendre chi function. We use  equation (\ref{eq:second_int}) and set $a=0,b = z^2, k = s, m = 1/2$ and simplify using equation (5) in \cite{guillera}.
\begin{equation}
\int_0^1 \frac{\log ^{-1+s}\left(\frac{1}{x}\right)}{\sqrt{x} \left(1-x z^2\right)} \, dx=\Gamma (s) \Phi
   \left(z^2,s,\frac{1}{2}\right)
\end{equation}
where $s,z\in\mathbb{C}$.
\end{example}
\begin{example}
In this example we use  equation (\ref{eq:second_int}) and replace $x\to t^a,m\to z/a$ and simplify.
\begin{equation}\label{eq:6.29}
\int_0^1 \frac{x^{z-1} \log ^{k-1}\left(\frac{1}{x}\right)}{1-b x^a} \,
   dx=a^{-k} \Gamma (k) \Phi \left(b,k,\frac{z}{a}\right)
\end{equation}
where $Re(b)<0$.
\end{example}
\begin{example}
In this example we use  equation (\ref{eq:6.29}) and replace $z\to z+1$, then form a second equation by replacing $z\to s$ and take their difference. Next set $b=-1$ and simplify in terms of the Hurwitz zeta using equation (64:13:3) in \cite{atlas}. Next we take the first partial derivative with respect to $k$ and apply l'Hopital's rule as $k\to 0$ and simplify in terms of the log-gamma function, [DLMF,\href{https://dlmf.nist.gov/25.11.E18}{25.11.18}] and Euler's constant [DLMF,\href{https://dlmf.nist.gov/5.4.E11}{5.4.11}].
\begin{multline}
\int_0^1 \frac{\left(x^s-x^z\right) \log \left(\log \left(\frac{1}{x}\right)\right)}{\left(1+x^a\right) \log
   \left(\frac{1}{x}\right)} \, dx
=\frac{1}{2} \left(\log (4) \log \left(\Gamma \left(\frac{1-2 a+z}{2
   a}\right)\right)-2 \log (2 a \exp (\gamma )) \right. \\ \left.
\log \left(\frac{2 a \Gamma \left(\frac{1+s}{2 a}\right) \Gamma
   \left(\frac{1+a+z}{2 a}\right)}{(1-2 a+z) \Gamma \left(\frac{1+a+s}{2 a}\right)}\right)+2 (\gamma +\log (a)) \log
   \left(\Gamma \left(-1+\frac{1+z}{2 a}\right)\right)\right. \\ \left.
+\zeta''\left(0,\frac{1+s}{2
   a}\right)-\zeta''\left(0,\frac{1+a+s}{2 a}\right)
-\zeta''\left(0,\frac{1+z}{2a}\right)+\zeta''\left(0,\frac{1+a+z}{2 a}\right)\right)
\end{multline}
where $0< Re(s)<1,0< Re(z)<1$.
\end{example}
\begin{example}
In this example we use  equation (\ref{eq:6.29}) and replace $z\to z+1$, then form a second equation by replacing $z\to s$ and take their difference. Next set $b=-1$ and simplify in terms of the Hurwitz zeta using equation (64:13:3) in \cite{atlas}. Next apply l'Hopital's rule as $k\to 0$ and simplify in terms of the log-gamma function using [DLMF,\href{https://dlmf.nist.gov/25.11.E18}{25.11.18}].
\begin{equation}
\int_0^1 \frac{x^m-x^s}{(1+x) \log \left(\frac{1}{x}\right)} \, dx=\log \left(\frac{\Gamma
   \left(\frac{1+m}{2}\right) \Gamma \left(1+\frac{s}{2}\right)}{\Gamma \left(1+\frac{m}{2}\right) \Gamma
   \left(\frac{1+s}{2}\right)}\right)
\end{equation}
where $0< Re(m)<1,0< Re(s)<1$
\end{example}
\begin{example}
Derivation of a generalized form for equation (2.6.18.6) in \cite{prud1} in terms of the incomplete gamma function.  We use equation (\ref{eq:6.5}) and form a second equation by replacing $b\to -b$ and taking their difference. Next replace $x\to t^u$, set $m=1/2$ and simplify.
\begin{multline}\label{eq:6.32}
\int_0^1 \frac{x^{-1+\frac{u}{2}}}{\left(-1+b x^u\right) \left(a^2+\log ^2(x)\right)} \, dx\\
=\sum _{n=0}^{\infty
   } \frac{i b^n e^{-i a \left(\frac{1}{2}+n\right) u} \left(\Gamma \left(0,-i a \left(\frac{1}{2}+n\right)
   u\right)-e^{i a (1+2 n) u} \Gamma \left(0,i a \left(\frac{1}{2}+n\right) u\right)\right)}{2 a}
\end{multline}
where $Re(a)>0,Re(u)>0$.
\end{example}
\begin{example}
Derivation of equation (2.6.18.6) in \cite{prud1} in terms of the incomplete gamma function. We use equation (\ref{eq:6.32}) and set $b=-1$.
\begin{multline}
\int_0^1 \frac{x^{-1+\frac{u}{2}}}{\left(1+x^u\right) \left(a^2+\log ^2(x)\right)} \, dx\\
=\sum _{n=0}^{\infty }
   \frac{i (-1)^n e^{-i a \left(\frac{1}{2}+n\right) u} \left(-\Gamma \left(0,-i a \left(\frac{1}{2}+n\right)
   u\right)+e^{i a (1+2 n) u} \Gamma \left(0,i a \left(\frac{1}{2}+n\right) u\right)\right)}{2 a}
\end{multline}
where $Re(a)>0,Re(u)>0$.
\end{example}
\begin{example}
Derivation of the digamma function in terms of the infinite sum involving the incomplete gamma function. Since the left-hand sides of equations (\ref{eq:6.32}) and (2.6.18.6) in \cite{prud1}, we can equate the right-hand sides to yield the stated result.
\begin{multline}
\sum _{n=0}^{\infty } i^{2 n+1} \left(e^{i a \left(\frac{1}{2}+n\right) u} \Gamma \left(0,i a
   \left(\frac{1}{2}+n\right) u\right)-e^{-i a \left(\frac{1}{2}+n\right) u} \Gamma \left(0,-i a
   \left(\frac{1}{2}+n\right) u\right)\right)\\
=\frac{1}{2} \left(-\psi ^{(0)}\left(\frac{\pi +a u}{4 \pi }\right)+\psi
   ^{(0)}\left(\frac{1}{2} \left(1+\frac{\pi +a u}{2 \pi }\right)\right)\right)
\end{multline}
where $Re(a)>0,Re(u)>0$.
\end{example}
\begin{example}
Derivation of a generalized form of equation (2.6.18.11-14) in \cite{prud1}. We use equation (\ref{eq:second_int}) and set $m=m+1,k=-1,a=e^{a}$, next we form a second equation by replacing $a\to -a$ and add these two equations. Next we replace $t \to -i \log(b), a \to a i$, and using this equation we replace $b\to -b$ and take their difference and simplify.
\begin{multline}\label{eq:6.35}
\int_0^1 \frac{x^{1+m} \log (x)}{\left(1-b^2 x^2\right) \left(a^2+\log ^2(x)\right)} \, dx\\
=\sum _{n=0}^{\infty
   } \frac{\left((-b)^n-b^n\right) e^{-i a (1+m+n)} \left(\Gamma (0,-i a (1+m+n))+e^{2 i a (1+m+n)} \Gamma (0,i a (1+m+n))\right)}{4 b}
\end{multline}
where $-1/2\leq Re(b)\leq 1/2$.
\end{example}
\begin{example}
Derivation Derivation of a generalized form of equation (2.6.18.11-14) in \cite{prud1}. In this example we simply set $m=-1$ in equation (\ref{eq:6.35}) and simplify.
\begin{multline}
\int_0^1 \frac{\log (x)}{\left(1-b^2 x^2\right) \left(a^2+\log
   ^2(x)\right)} \, dx\\
=\sum _{n=1}^{\infty } \frac{1}{4} \left(-1+(-1)^n\right)
   b^{-1+n} \left(e^{-i a n} \Gamma (0,-i a n)+e^{i a n} \Gamma (0,i a
   n)\right)
\end{multline}
where $-1/2\leq Re(b)\leq 1/2$.
\end{example}
\begin{example}
 We use equation (\ref{eq:second_int}) and set $m=m+1,k=-1,a=e^{a}$, next we form a second equation by replacing $a\to -a$ and add these two equations.Next we take the $j$-th derivative with respect to $a$ and simplify.
\begin{multline}
\int_0^1 \frac{x^m \left((-a-i \log (x))^{-1-j}-(-a+i \log (x))^{-1-j}\right)}{-1+b x} \, dx\\
=-\sum _{n=0}^{\infty
   } i b^n e^{-i a (1+m+n)} \left(\left(-\frac{i}{a (1+m+n)}\right)^j (-i (1+m+n))^j E_{1+j}(-i a (1+m+n))\right. \\ \left.
+e^{2 i a
   (1+m+n)} \left(\frac{i}{a (1+m+n)}\right)^j (i (1+m+n))^j E_{1+j}(i a (1+m+n))\right)
\end{multline}
where $Re(a)>0,-1/2\leq Re(b)\leq 1/2$.
\end{example}
\begin{example}
In this example we use equation (\ref{eq:second_int}) and replace $t\to -i\log(b),a\to e^a$, next we set $k=-1$ and take the $p$-th derivative with respect to $b$ and the $q$-th derivative with respect $a$ and simplify using [Wolfram,\href{http://functions.wolfram.com/06.06.20.0013.01}{01}].
\begin{multline}
\int_0^1 \frac{x^{-1+m+p}}{(1-b x)^{p+1} (a+\log (x))^{q+1}} \, dx\\
=\sum _{n=0}^{\infty } \frac{(-1)^{-2 p-q+1}
   b^{n-p} e^{-a (m+n)} (-m-n)^q (a (m+n))^{-q} E_{1+q}(-a (m+n)) \Gamma (1+n)}{\Gamma (1+n-p) \Gamma (1+p)}
\end{multline}
where $Re(a)>0,-1/2\leq Re(b)\leq 1/2,p,q\in\mathbb{Z_{+}}$.
\end{example}
\begin{example}
In this example we use equation (\ref{eq:second_int}) and replace $t\to -i\log(b),a\to e^a$, next we take the first partial derivative with respect to $k$ and set $k=-1$. Next we take the $p$-th derivative with respect to $b$ and simplify using [Wolfram,\href{http://functions.wolfram.com/06.06.20.0013.01}{01}].
\begin{multline}
\int_0^1 \frac{x^{-1+m} \log (a+\log (x))}{(a+\log (x)) (1-b x)^{p+1}} \, dx\\
=\sum _{n=0}^{\infty } \frac{b^{n-p}
    \Gamma (1+n) \left(\Gamma (0,-a (m+n-p)) \log \left(\frac{1}{a}\right)-G_{2,3}^{3,0}\left(-a
   (m+n-p)\left|
\begin{array}{c}
 1,1 \\
 0,0,0 \\
\end{array}
\right.\right)\right)}{e^{a (m+n-p)}\Gamma (1+n-p) \Gamma (1+p)}
\end{multline}
where $Re(a)>0,-1/2\leq Re(b)\leq 1/2,p\in\mathbb{Z_{+}}$.
\end{example}
\begin{example}
In this example we use equation (\ref{eq:second_int}) and replace $t\to -i\log(b)$ and set $a=0$. Next we replace $x\to x^q,b\to -b,m\to z/q$ and simplify and take the $j$-th derivative with respect to $b$ and simplify the factorial factors.
\begin{multline}
\int_0^1 \frac{(-\log (x))^{-1+r} x^{-1+z}}{\left(1+b x^q\right)^j} \, dx
=\sum _{n=0}^{\infty }
   \frac{(-b)^{1-j+n} q^{-r} \left(1-j+n+\frac{z}{q}\right)^{-r} \Gamma (1+n) \Gamma (r)}{\Gamma (j) \Gamma
   (2-j+n)}
\end{multline}
where $Re(z)>0,|Re(b)|<1,j\in\mathbb{Z_{+}}$.
\end{example}
\begin{example}
Extended form of section (2.6.19) in \cite{prud1}. In this example we use equation (\ref{eq:second_int}) and replace $t\to -i\log(b)$ and set $a=0$. Next we replace $x\to x^q,b\to-b,m\to z/q$. Next we replace $z\to z+q$ and take the indefinite integral with respect to $b$.
\begin{multline}\label{eq:6.41}
\int_0^1 x^{z-1} \log ^{k-1}(x) \log \left(1+b x^q\right) \, dx=(-1)^{k-1} b \Gamma (k) \sum _{n=0}^{\infty }
   \frac{(-b)^n}{(n+1) ((1+n) q+z)^k}
\end{multline}
where $Re(q)>0,Re(k)>0,Re(z)>0,|Re(b)|<1$.
\end{example}
\begin{example}
Derivation of equation (2.6.19.6) in \cite{prud1}.  Use equation (\ref{eq:6.41}) and set $z=0,k\to k+1$ and simplify the sum in terms of the Polylogarithm function using equation [DLMF,\href{https://dlmf.nist.gov/25.12.E1}{25.12.1}].
\begin{equation}
\int_0^1 \frac{\log ^k\left(\frac{1}{x}\right) \log \left(1+b x^q\right)}{x} \, dx=-q^{-1-k} \Gamma (1+k)
   \text{Li}_{2+k}(-b)
\end{equation}
where $Re(k)>0,Re(q)>0,|Re(b)|<1$.
\end{example}
\begin{example}
Generalized equation (2.6.19.7) in \cite{prud1}. Use equation (\ref{eq:6.41}) and replace $q\to u$. We form two equations by replacing $b\to \frac{1}{\sqrt{b^2 \left(c^2-1\right)}+b c}$ and $b\to \frac{\sqrt{b^2 \left(c^2-1\right)}+b c}{b^2}$ and add these two equations and simplify.
\begin{multline}\label{eq:6.43}
\int_0^1 x^{-1+z} \log ^{-1+k}(x) \log \left(1+b^2 x^{2 u}+2 b x^u \cos (t)\right) \, dx\\
=2 \Gamma (k) (-1)^k
   \sum _{n=1}^{\infty } \frac{(-b)^n \cos (n t)}{n (n u+z)^k}
\end{multline}
where $Re(k)>0,Re(u)>0,Re(z)>0,|Re(b)|<1$.
\end{example}
\begin{example}
Here we use equation (\ref{eq:6.43}) and replace $n\to n+1$ and form a second equation by replacing $b\to -b$ and take their difference and simplify the left-hand side using equation [Wolfram,\href{https://mathworld.wolfram.com/InverseTangent.html}{1}].
\begin{equation}\label{eq:6.44}
\int_0^1 x^{z-1} \tanh ^{-1}\left(b x^u\right) \log ^{k-1}\left(\frac{1}{x}\right) \, dx=\frac{1}{2} \Gamma (k)
   \sum _{n=1}^{\infty } \frac{\left(1-(-1)^n\right) b^n}{n (n u+z)^k}
\end{equation}
where $Re(k)>0,Re(u)>0,|Re(b)|<1$.
\end{example}
\begin{example}
Here we use equation (\ref{eq:6.44}) and set $z=0$ and simplify the right-hand side using equation [DLMF,\href{https://dlmf.nist.gov/25.12.E1}{25.12.1}].
\begin{equation}\label{eq:6.45}
\int_0^1 \frac{\tanh ^{-1}\left(b x^u\right) \log ^{-1+k}\left(\frac{1}{x}\right)}{x} \, dx=-\frac{1}{2} u^{-k}
   \Gamma (k) (\text{Li}_{1+k}(-b)-\text{Li}_{1+k}(b))
\end{equation}
where $Re(k)>0,Re(u)>0,|Re(b)|<1$.
\end{example}
\begin{example}
Here we use equation (\ref{eq:6.45}) and simplify using equation [DLMF,\href{https://dlmf.nist.gov/25.12.E13}{25.12.13}].
\begin{multline}
\int_0^1 \frac{\tanh ^{-1}\left(a x^u\right) \log ^{-2+s}\left(\frac{1}{x}\right)}{x} \, dx\\
=\frac{i^{1-s}
   2^{-2+s} \pi ^s u^{1-s} \csc (\pi  s) }{-1+s}\left(\zeta \left(1-s,\frac{\pi -i \log (-a)}{2 \pi }\right)-\zeta
   \left(1-s,\frac{\pi -i \log (a)}{2 \pi }\right)\right. \\ \left.
+i^{2 s} \left(-\zeta \left(1-s,\frac{\pi +i \log (-a)}{2 \pi
   }\right)+\zeta \left(1-s,\frac{\pi +i \log (a)}{2 \pi }\right)\right)\right)
\end{multline}
where $Re(s)>1,|Re(u)|<1,Re(a)>1$.
\end{example}
\begin{example}
Extended form of equation (2.6.18.16) in \cite{prud1}. In this example we use equation (\ref{eq:second_int}) and replace $t\to -i\log(b)$ and replace $x\to x^u,a\to a u,m\to s/u$. Next we form a second equation by replacing $a\to -a$ and take their difference and set $k=-1,a=a i$. Next take the indefinite integral with respect to $b$ and replace $s\to s+u$ and simplify.
   \begin{multline}
\int_0^1 \frac{x^{-1+s} \log \left(1-b x^u\right)}{a^2+\log ^2(x)} \, dx\\
=\sum _{n=0}^{\infty } \frac{i b^{1+n}
   \left(e^{-i a (s+u+n u)} \Gamma (0,-i a (s+u+n u))-e^{i a (s+u+n u)} \Gamma (0,i a (s+u+n u))\right)}{2 a
   (1+n)}
\end{multline}
where $Re(a)>0,Re(u)>0,Re(s)>0,|Re(b)|<1$.
\end{example}
\section{Extended Prudnikov integral forms}
In this section we look at deriving formulae using the definite integral substitution. This method is very useful in deriving finite series, products and functional equations. This is a continuation of the work done in the previous section.
When we apply the contour integral method in \cite{reyn4} to equation (2.2.9.1) in \cite{prud1} we get the following theorem;
\begin{theorem}
For all $Re(w)>0,Re(k)>0,Re(w)>0,Re(m)>0,|\gamma|<\pi$,
\begin{multline}\label{eq:prud1}
\frac{1}{2\pi i}\int_{C}\int_{0}^{1}\frac{a^w w^{-k-1} x^{m+w-1}}{x^2+2 x \cos (\gamma )+1}dxdw\\
=\frac{1}{2\pi i}\int_{C}\sum_{j=0}^{\infty}\frac{(-1)^j a^w w^{-k-1} (\cos (\gamma  j)+\cot (\gamma
   ) \sin (\gamma  j))}{j+m+w}dw
\end{multline}
\end{theorem}
\subsection{Derivation of the left-hand side contour integral representation}
Using a generalization of Cauchy's integral formula \ref{intro:cauchy}, we form the definite integral by replacing $y$ by $\log{ax}$ and multiply both sides by $\frac{x^{m-1}}{x^2+2 x \cos (\gamma )+1}$ to get;
\begin{multline}\label{eq:prud2}
\int_{0}^{1}\frac{x^{m-1} \log ^k(a x)}{k! \left(x^2+2 x \cos (\gamma )+1\right)}dx\\
=\frac{1}{2\pi i}\int_{0}^{1}\int_{C}\frac{w^{-k-1} x^{m-1} (a x)^w}{x^2+2 x
   \cos (\gamma )+1}dwdx\\
   =\frac{1}{2\pi i}\int_{C}\int_{0}^{1}\frac{w^{-k-1} x^{m-1} (a x)^w}{x^2+2 x
   \cos (\gamma )+1}dxdw
\end{multline}
We are able to switch the order of integration over $x$ and $w$ using Fubini's theorem for multiple integrals see page 178 in \cite{gelca}, since the integrand is of bounded measure over the space $\mathbb{C} \times [0,1]$.
\subsection{Derivation of the right-hand side contour integral representation}
In this section, we will once again use Cauchy's generalized integral formula, equation (\ref{intro:cauchy}), and take the infinite integral to derive equivalent series representation for the contour integrals. We proceed using equation~(\ref{intro:cauchy}) and replace $y$ by $\log (a)+x$ and multiply both sides by $e^{m x}$ and simplify. Next, we take the definite integral over 
$x\in]0,\infty)$ and simplify using equation (3.382.4) in \cite{grad}. Next we replace $m\to j+m$ and multiply both sides by $-(-1)^j (\cos (\gamma  j)+\cot (\gamma ) \sin (\gamma  j))$ then take the infinite sum over $j\in[0,\infty)$ and simplify in terms of the incomplete gamma function to obtain
\begin{multline}\label{eq:prud3}
-\sum_{j=0}^\infty{}\frac{(-1)^j a^{-j-m} (-j-m)^{-1-k} \Gamma (1+k,-((j+m) \log (a))) (\cos (j \gamma )+\cot (\gamma ) \sin (j
   \gamma ))}{k!}\\
   =\frac{1}{2\pi i}\sum_{j=0}^{\infty}\int_{C}\frac{(-1)^j a^w w^{-1-k} (\cos (j \gamma )+\cot (\gamma ) \sin (j \gamma ))}{j+m+w}dw\\
=\frac{1}{2\pi i}\int_{C}\sum_{j=0}^{\infty}\frac{(-1)^j a^w w^{-1-k} (\cos (j \gamma )+\cot (\gamma ) \sin (j \gamma ))}{j+m+w}dw\\
\end{multline}
We are able to switch the order of integration and summation over $w$ using Tonellii's theorem for  integrals and sums see page 177 in \cite{gelca}, since the summand is of bounded measure over the space $\mathbb{C} \times [0,\infty)$.
\begin{theorem}
For all $Re(m)>0,|\gamma|<\pi$ then,
\begin{multline}\label{eq:prud_int}
\int_0^1 \frac{x^{-1+m} \log ^k(a x)}{1+x^2+2 x \cos (\gamma )} \, dx\\
=-\sum _{j=0}^{\infty } (-1)^j a^{-j-m}
   (-j-m)^{-1-k} \Gamma (1+k,-((j+m) \log (a))) (\cos (j \gamma )+\cot (\gamma ) \sin (j \gamma ))
\end{multline}
\end{theorem}
\begin{proof}
Since the right-hand sides of equations (\ref{eq:prud2}) and (\ref{eq:prud3}) are equivalent relative to equation (\ref{eq:prud1}), we can equate the left-hand sides and simplify the gamma function to yield the stated result.
\end{proof}
\begin{example}
In this example we use equation (\ref{eq:prud_int}) and take the indefinite integral with respect to $b$ and simplify;
\begin{multline}\label{eq:7.5}
\int_0^1 x^{-2+m} \log ^k(a x) \log \left(1+2 b x+x^2\right) \, dx\\
=-\sum _{j=0}^{\infty } \frac{2 (-1)^j
    \Gamma (1+k,-((j+m) \log (a))) \left(b \cos \left(j \cos ^{-1}(b)\right)-\sqrt{1-b^2} \sin \left(j \cos
   ^{-1}(b)\right)\right)}{a^{j+m}(1+j) (-1)^{k+1} (j+m)^{k+1}}
\end{multline}
where $Re(a)>0,0< Re(b)<1,Re(m)>0$.
\end{example}
\begin{example}
In this example we use equation (\ref{eq:prud_int}) set $a=1$ and take the indefinite integral with respect to $b$ and simplify the gamma function using [Wolfram,\href{https://mathworld.wolfram.com/IncompleteGammaFunction.html}{1}];
\begin{multline}\label{eq:7.6}
\int_0^1 x^{-2+m} \log ^k\left(\frac{1}{x}\right) \log \left(1+2 b x+x^2\right) \, dx\\
=2 \Gamma (1+k) \sum
   _{j=0}^{\infty } \frac{(-1)^j \left(b \cos \left(j \cos ^{-1}(b)\right)-\sqrt{1-b^2} \sin \left(j \cos
   ^{-1}(b)\right)\right)}{(1+j) (j+m)^{k+1}}
\end{multline}
where $Re(m)>0,0< Re(b)<1$.
\end{example}
\begin{example}
In this example we use equation (\ref{eq:7.6}) and set $k=0,b=0,m\to m+1$. Next replace $x\to x^p$ and simplify. Then replace $m\to q/p-1$ followed by $p\to p/2$. The final integral form is equivalent to equation 325(2) in \cite{grobner} and hence equating the right-hand sides yields the stated result. The functions used are the digamma function [DLMF,\href{https://dlmf.nist.gov/5.15.E5}{5.15.1}] and  hypergeometric function [DLMF,\href{https://dlmf.nist.gov/15.2}{15.2}].
\begin{multline}
\psi ^{(0)}\left(\frac{p+q}{2 p}\right)-\psi ^{(0)}\left(1+\frac{q}{2 p}\right)\\
=-\frac{2 i p \left(\,
   _2F_1\left(1,1+\frac{2 q}{p};\frac{2 (p+q)}{p};-i\right)-\, _2F_1\left(1,1+\frac{2 q}{p};\frac{2
   (p+q)}{p};i\right)\right)}{p+2 q}
\end{multline}
where $Re(q)>0,Re(p)>0$.
\end{example}
\begin{example}
In this example we use equation (\ref{eq:prud_int}) and set $k\to 2k,b=0$. Next replace $x\to x^p$ and simplify. Next replace $m\to z/p,p\to p/2$ and simplify.
\begin{multline}\label{eq:7.8}
\int_0^1 x^{-1+z} \log ^{2 k}\left(\frac{1}{x}\right) \log \left(1+x^p\right) \, dx\\
=-\sum _{j=0}^{\infty }
   \frac{(-1)^j 2^{2+2 k} p^{-1-2 k} \left(1+j+\frac{2 z}{p}\right)^{-1-2 k} \Gamma (1+2 k) \sin \left(\frac{j \pi
   }{2}\right)}{1+j}
\end{multline}
where $Re(z)>0$.
\end{example}
\begin{example}
In this example we derive a definite integral in terms of Bernoulli numbers also studied by Gr\"{o}ber  see equation 326(3) in \cite{grobner}. We use equation (\ref{eq:7.8}) and set $z=0$ and simplify the series in terms of the  polylogarithm function using [DLMF,\href{https://dlmf.nist.gov/25.13.E3}{25.13.2}], and the Hurwitz-zeta function in terms of Bernoulli numbers using  [Wolfram,\href{https://mathworld.wolfram.com/RiemannZetaFunction.html}{66}].
\begin{equation}
\int_0^1 \frac{\log ^{2 k}\left(\frac{1}{x}\right) \log \left(1+x^p\right)}{x} \, dx=\frac{(-1)^k \left(2^{2
   k+1}-1\right) \pi ^{2+2 k} B_{2+2 k}}{2 p^{2 k+1} (1+k) (1+2 k)}
\end{equation}
where $Re(p)>0$.
\end{example}
\begin{example}
In this example we use equation (\ref{eq:7.5}) and set $k=0,a=1,b\to \cos(b),m=1$ and simplify.
\begin{multline}
\int_0^1 \log \left(1+x^2+2 x \cos (b)\right) \, dx=(1+\cos (b)) \log (2 (1+\cos (b)))+b \sin (b)-2
\end{multline}
where $Re(b)>0$.
\end{example}
\begin{example}
In this example we use equation (\ref{eq:7.5}) set $k=-1,a\to e^a,m\to m+1$. Next form a second equation by replacing $a\to -a$ and take their difference.
\begin{multline}
\int_0^1 \frac{x^{-1+m} \log \left(1+2 b x+x^2\right)}{a^2+\log ^2(x)} \, dx\\
=\sum _{j=0}^{\infty } \frac{i
   (-1)^j }{a (1+j)}\left(e^{i a (1+j+m)} \Gamma (0,i a (1+j+m))-e^{-i a (1+j+m)} \Gamma (0,-i a (1+j+m))\right) \\
\left(b \cos
   \left(j \cos ^{-1}(b)\right)-\sqrt{1-b^2} \sin \left(j \cos ^{-1}(b)\right)\right)
\end{multline}
where $Re(a)>0,0< Re(b)<1,Re(m)>0$.
\end{example}
\begin{example}
Derivation of a functional equation involving the incomplete Beta function. In this example we use equation (2.6.14.18) in \cite{prud1} and rewrite the infinite integral over [0,1] and simplify. We start with equation (\ref{eq:7.5}) with $a=1,k=0$. This gives an integral in terms of the incomplete Beta function from equation [DLMF,\href{https://dlmf.nist.gov/8.17.E1}{8.17.1}]. Next using equation (2.6.14.18) in \cite{prud1}, we transform the integral over $[0,\infty) \to [0,1]$ and replace with the previously derived integral and simplify to yield the stated result.
\begin{multline}
\{b^{2 \alpha } \left(B_{\frac{1}{b}}(1+\alpha ,0)-B_b(1-\alpha ,0)\right)+B_b(1+\alpha
   ,0)-B_{\frac{1}{b}}(1-\alpha ,0)\\
=i \left(-1+b^{2 \alpha }\right) \pi -\frac{2 b^{\alpha }}{\alpha }+\left(1+b^{2
   \alpha }\right) \pi  \cot (\pi  \alpha )
\end{multline}
where $Re(b)>0,Re(\alpha)>0$.
\end{example}
\begin{example}
Derivation of equation (2.6.17.15) in \cite{prud1} for complex numbers. In this example we use equation (\ref{eq:prud_int}) and set $a=1,m\to m+1$ and simplify in terms of the Hurwitz-Lerch zeta function with $\gamma\to b $. Next we form three equations when $m=2,m=1,m=0$ and add all three equations and simplify. Next simplify the polylogarithm function using [DLMF,\href{https://dlmf.nist.gov/25.13.E3}{25.13.2}]. Next we apply l'Hopital's rule as $k\to -1$ and simplify. Next we simplify the first partial derivative of the Huritz-Lerch zeta function using equation (\ref{eq:fun_eq}) and simplify to get the stated result. The functions used in this example are Euler's constant $\gamma$ in [Wolfram,\href{http://functions.wolfram.com/02.06.02.0001.01}{1}], Stieltjes constant $\gamma_{n}$ in [Wolfram,\href{http://functions.wolfram.com/10.05.02.0001.01}{1}] and the digamma function $\psi(z)$ in [Wolfram,\href{http://functions.wolfram.com/06.14.02.0001.01}{1}].
\begin{multline}\label{eq:7.13}
\int_0^1 \frac{(-1+x)^2}{\left(1+x^2+2 x \cos (b \pi )\right) \log (x)} \, dx\\
=\frac{2 e^{\frac{i b \pi }{2}}}{\left(-1+e^{i b \pi }\right) \left(1+e^{i b \pi }\right)^2 \pi }
   \left(-i e^{i b \pi } \pi  (\gamma +\log (2 \pi )) \sin \left(\frac{b \pi }{2}\right)\right. \\ \left.
+\cos \left(\frac{b \pi
   }{2}\right) \left(-i e^{i b \pi } \gamma _1\left(\frac{1-b}{2}\right)+i e^{i b \pi } \gamma
   _1\left(\frac{1+b}{2}\right)-\left(2+e^{i b \pi }\right) \pi  \Phi'\left(-e^{-i b \pi
   },0,2\right)\right.\right. \\ \left.\left.
+\left(e^{i b \pi }+2 e^{2 i b \pi }\right) \pi  \Phi'\left(-e^{i b \pi
   },0,2\right)\right)\right)\\
   =\gamma +\log (2 \pi )+\cos ^2\left(\frac{b \pi }{2}\right) \left(\psi ^{(0)}\left(\frac{1}{2}-\frac{b}{2}\right)+\psi ^{(0)}\left(\frac{1+b}{2}\right)\right)\\
+\frac{\cos (b \pi ) \cot
   \left(\frac{b \pi }{2}\right) \left(-\gamma _1\left(\frac{1-b}{2}\right)+\gamma _1\left(\frac{1+b}{2}\right)\right)}{\pi }
\end{multline}
where $Re(b)>0$.
\end{example}
\begin{example}
In this example we write equation (\ref{eq:7.13}) in terms of integers and compare and equate to equation (2.6.17.15) in \cite{prud1}. Where $b\to m/n$.
\begin{multline}\label{eq:7.14}
\int_0^1 \frac{(-1+x)^2}{\left(1+x^2+2 x \cos \left(\frac{m \pi }{n}\right)\right) \log (x)} \, dx\\
=\sum
   _{k=1}^{\frac{n-1}{r}}\frac{ (-1)^k \sin \left(\frac{k m \pi }{n}\right) \log \left(\frac{\Gamma \left(\frac{r
   \left(n-(-1)^r k+1\right)}{2 n}\right) \Gamma \left(\frac{r \left(n-(-1)^r k+1\right)}{2 n}\right) \Gamma
   \left(\frac{r (k+2)}{2 n}\right) \Gamma \left(\frac{r k}{2 n}\right)}{\Gamma \left(\frac{r (k+1)}{2 n}\right) \Gamma
   \left(\frac{r (k+1)}{2 n}\right) \Gamma \left(\frac{r \left(n-(-1)^r k\right)}{2 n}\right) \Gamma \left(\frac{r
   \left(n-(-1)^r k+2\right)}{2 n}\right)}\right)}{\sin \left(\frac{m \pi }{n}\right)}\\
=\gamma +\log (2 \pi )+\cos
   ^2\left(\frac{m \pi }{2 n}\right) \left(\psi ^{(0)}\left(\frac{1}{2}-\frac{m}{2 n}\right)+\psi
   ^{(0)}\left(\frac{1}{2} \left(1+\frac{m}{n}\right)\right)\right)\\+\frac{\cos \left(\frac{m \pi }{n}\right) \cot
   \left(\frac{m \pi }{2 n}\right) \left(-\gamma _1\left(\frac{1}{2} \left(1-\frac{m}{n}\right)\right)+\gamma
   _1\left(\frac{1}{2} \left(1+\frac{m}{n}\right)\right)\right)}{\pi }
\end{multline}
where $m<n$.
\end{example}
\begin{example}
In this example we use equation (\ref{eq:7.14}) and equate the right-hand sides and simplify to yield the stated result. Note this is a special case of equation (2.6.17.15) in \cite{prud1} when $r=1$.
\begin{multline}\label{eq:7.15}
\sum _{k=1}^{n-1} (-1)^k \log \left(\frac{\Gamma \left(\frac{k}{2 n}\right) \Gamma \left(\frac{2+k}{2 n}\right) \Gamma \left(\frac{1+k+n}{2 n}\right)^2}{\Gamma \left(\frac{1+k}{2 n}\right)^2 \Gamma \left(\frac{k+n}{2 n}\right) \Gamma \left(\frac{2+k+n}{2
   n}\right)}\right) \sin \left(\frac{k m \pi }{n}\right)\\
=\sin \left(\frac{m \pi }{n}\right) \left(\gamma +\log (2 \pi )+\cos ^2\left(\frac{m \pi }{2 n}\right) \left(\psi ^{(0)}\left(\frac{-m+n}{2 n}\right)+\psi ^{(0)}\left(\frac{m+n}{2 n}\right)\right)\right. \\ \left.
+\frac{\cos
   \left(\frac{m \pi }{n}\right) \cot \left(\frac{m \pi }{2 n}\right) \left(-\gamma _1\left(\frac{-m+n}{2 n}\right)+\gamma _1\left(\frac{m+n}{2 n}\right)\right)}{\pi }\right)
\end{multline}
where $|Re(m/n)|<\pi/4$.
\end{example}
\begin{example}
A finite product involving the Gamma function in terms of fundamental constants. In this example we use equation (\ref{eq:7.15}) and take the exponential function of both sides and simplify.
\begin{multline}
\prod _{k=1}^{n-1} \left(\frac{\Gamma \left(\frac{k}{2 n}\right) \Gamma \left(\frac{2+k}{2 n}\right) \Gamma \left(\frac{1+k+n}{2 n}\right)^2}{\Gamma \left(\frac{1+k}{2 n}\right)^2 \Gamma \left(\frac{k+n}{2 n}\right) \Gamma \left(\frac{2+k+n}{2 n}\right)}\right)^{(-1)^k
   \sin \left(\frac{k m \pi }{n}\right)}\\
=\exp \left(\sin \left(\frac{m \pi }{n}\right) \left(\gamma +\log (2 \pi )+\cos ^2\left(\frac{m \pi }{2 n}\right) \left(\psi ^{(0)}\left(\frac{-m+n}{2 n}\right)+\psi ^{(0)}\left(\frac{m+n}{2 n}\right)\right)\right.\right. \\ \left.\left.
+\frac{\cos \left(\frac{m \pi
   }{n}\right) \cot \left(\frac{m \pi }{2 n}\right) \left(-\gamma _1\left(\frac{-m+n}{2 n}\right)+\gamma _1\left(\frac{m+n}{2 n}\right)\right)}{\pi }\right)\right)
\end{multline}
where $|Re(m/n)|<\pi/4$.
\end{example}
\begin{example}
Extended form of equations (2.6.7.23-24) in \cite{prud1}. In this example we use equation (\ref{eq:7.6}) and replace $m\to m+1$. Next replace $m \to s$ to form a second equation and take their difference. Next apply l'Hopital's as $k\to -1$ and simplify.
\begin{multline}\label{eq:7.17}
\int_0^1 \frac{\left(x^m-x^s\right) \log \left(1+2 b x+x^2\right)}{\log (x)} \, dx\\
=2 \sum _{j=0}^{\infty }
   \frac{(-1)^j \log \left(\frac{2+j+m}{2+j+s}\right) \left(b \cos \left(j \cos ^{-1}(b)\right)-\sqrt{1-b^2} \sin
   \left(j \cos ^{-1}(b)\right)\right)}{1+j}
\end{multline}
where $Re(m)>0,Re(n)>0,Re(b)>0$.
\end{example}
\begin{example}
In this example we use equation (\ref{eq:7.17}) and set $b=0$ and replace $x \to x^p$ and simplify.
\begin{multline}
\int_0^1 \frac{\left(x^{\alpha }-x^{\beta }\right) \log \left(1+x^p\right)}{\log (x)} \, dx=-\sum _{j=0}^{\infty
   } \frac{2 (-1)^j \log \left(\frac{2+p+j p+2 \alpha }{2+p+j p+2 \beta }\right) \sin \left(\frac{j \pi
   }{2}\right)}{1+j}
\end{multline}
where $Re(\alpha)>0,Re(\beta)>0,Re(p)>0$.
\end{example}
\begin{example}
In this example we use equation (\ref{eq:7.6}) replace $m\to m+1$ and form a second equation by replacing $m\to s$ and taking their difference. Next take the first partial derivative with respect to $k$ and apply l'Hopital's rule as $k\to -1$ and simplify.
\begin{multline}
\int_0^1 \frac{\left(x^m-x^s\right) \log \left(1+2 b x+x^2\right) \log \left(\log
   \left(\frac{1}{x}\right)\right)}{\log \left(\frac{1}{x}\right)} \, dx\\
=\sum _{j=0}^{\infty } \frac{1}{1+j}(-1)^j \log
   \left(\frac{2+j+s}{2+j+m}\right) (2 \gamma +\log (2+j+m)+\log (2+j+s))\\ \left(-b \cos \left(j \cos
   ^{-1}(b)\right)+\sqrt{1-b^2} \sin \left(j \cos ^{-1}(b)\right)\right)
\end{multline}
where $Re(m)>0,Re(s)>0,Re(b)>0$.
\end{example}
\begin{example}
In this example we use equation (\ref{eq:7.6}) then take the first partial derivative with respect to $b$ and simplify in terms of the Hurwitz-Lerch zeta function using [DLMF,\href{https://dlmf.nist.gov/25.14.i}{25.14.1}]. Next form a second equation by replacing $m\to s$ and take their difference. Next apply l'Hopital's rule as $k\to -1$ and simplify.
\begin{multline}\label{eq:7.20}
\int_0^1 \frac{-x^m+x^s}{\left(1+x^2+2 x \cos (b)\right) \log \left(\frac{1}{x}\right)} \,
   dx\\
=\frac{\Phi'\left(-e^{-i b},0,1+m\right)-\Phi'\left(-e^{-i
   b},0,1+s\right)+e^{2 i b} \left(-\Phi'\left(-e^{i
   b},0,1+m\right)+\Phi'\left(-e^{i b},0,1+s\right)\right)}{-1+e^{2 i b}}
\end{multline}
where $Re(m)>0,Re(s)>0,Re(b)>0$.
\end{example}
\begin{example}
In this example we use equation (\ref{eq:7.20}) and simplify using equation (\ref{eq:fun_eq}).
\begin{multline}\label{eq:7.21}
\int_0^1 \frac{-x^m+x^s}{\left(1+x^2+2 x \cos (b)\right) \log \left(\frac{1}{x}\right)} \, dx\\
=-\frac{i}{4
   \left(-1+e^{2 i b}\right) \pi } \left(e^{i
   (b+b m-m \pi )} \left(\Phi \left(e^{-2 i m \pi },1,\frac{-b+\pi }{2 \pi }\right) (2 \gamma -i \pi +\log (4)+2 \log
   (\pi ))\right.\right. \\ \left.\left.
-2 \Phi'\left(e^{-2 i m \pi },1,\frac{-b+\pi }{2 \pi }\right)+e^{2 i m \pi }
   \left(-\Phi \left(e^{2 i m \pi },1,\frac{b+\pi }{2 \pi }\right) (2 \gamma +i \pi +\log (4)+2 \log (\pi ))\right.\right.\right. \\ \left.\left.\left.
+2
   \Phi'\left(e^{2 i m \pi },1,\frac{b+\pi }{2 \pi }\right)\right)\right)+e^{2 i b}
   \left(-e^{-i (b+b m+m \pi )} \left(\Phi \left(e^{-2 i m \pi },1,\frac{b+\pi }{2 \pi }\right) (2 \gamma -i \pi +\log
   (4)+2 \log (\pi ))\right.\right.\right. \\ \left.\left.\left.
-2 \Phi'\left(e^{-2 i m \pi },1,\frac{b+\pi }{2 \pi }\right)+e^{2 i m
   \pi } \left(-\Phi \left(e^{2 i m \pi },1,\frac{-b+\pi }{2 \pi }\right) (2 \gamma +i \pi +\log (4)+2 \log (\pi ))\right.\right.\right.\right. \\ \left.\left.\left.\left.
+2\Phi'\left(e^{2 i m \pi },1,\frac{-b+\pi }{2 \pi }\right)\right)\right)+e^{-i (b+(b+\pi )
   s)} \left(\Phi \left(e^{-2 i \pi  s},1,\frac{b+\pi }{2 \pi }\right) (2 \gamma -i \pi +\log (4)+2 \log (\pi ))\right.\right.\right. \\ \left.\left.\left.
-2\Phi'\left(e^{-2 i \pi  s},1,\frac{b+\pi }{2 \pi }\right)+e^{2 i \pi  s} \left(-\Phi
   \left(e^{2 i \pi  s},1,\frac{-b+\pi }{2 \pi }\right) (2 \gamma +i \pi +\log (4)+2 \log (\pi ))\right.\right.\right.\right. \\ \left.\left. \left.\left.\left.
+2 \Phi'\left(e^{2 i \pi  s},1,\frac{-b+\pi }{2 \pi }\right)\right)\right)\right)-e^{i (b+b
   s-\pi  s)} \left(\Phi \left(e^{-2 i \pi  s},1,\frac{-b+\pi }{2 \pi }\right) (2 \gamma -i \pi +\log (4)+2 \log (\pi
   ))\right.\right.\right. \\ \left.\left.
-2 \Phi'\left(e^{-2 i \pi  s},1,\frac{-b+\pi }{2 \pi }\right)+e^{2 i \pi  s} \left(-\Phi
   \left(e^{2 i \pi  s},1,\frac{b+\pi }{2 \pi }\right) (2 \gamma +i \pi +\log (4)+2 \log (\pi ))\right.\right.\right. \\ \left.\left.\left.
+2
   \Phi'\left(e^{2 i \pi  s},1,\frac{b+\pi }{2 \pi }\right)\right)\right)\right)
\end{multline}
where $Re(m)>0,Re(s)>0,Re(b)>0$.
\end{example}
\subsection{Special cases }
In these examples we evaluated equation (\ref{eq:7.21}) for $m=1/2,s=-1/2,b=\pi/2,\pi/3,\pi/4,\pi/6,\pi/8$ respectively.
\begin{example}
\begin{equation}
\int_0^1 \frac{1-x}{\sqrt{x} \left(1+x^2\right) \log \left(\frac{1}{x}\right)} \, dx=\log \left(\cot
   \left(\frac{\pi }{8}\right)\right)
\end{equation}
\end{example}
\begin{example}
\begin{equation}
\int_0^1 \frac{1-x}{\sqrt{x} \left(1+x+x^2\right) \log \left(\frac{1}{x}\right)} \, dx=\log (2)
\end{equation}
\end{example}
\begin{example}
\begin{multline}
\int_0^1 \frac{1-x}{\sqrt{x} \left(1+x \left(\sqrt{2}+x\right)\right) \log \left(\frac{1}{x}\right)} \, dx\\
=(1+i)
   (-1)^{5/8} \left(1+\sqrt[4]{-1}\right) \left(\cos \left(\frac{\pi }{8}\right) \log \left(\cot \left(\frac{3 \pi
   }{16}\right)\right)+\log \left(\tan \left(\frac{\pi }{16}\right)\right) \sin \left(\frac{\pi
   }{8}\right)\right)
\end{multline}
\end{example}
\begin{example}
\begin{multline}
\int_0^1 \frac{1-x}{\sqrt{x} \left(1+x \left(\sqrt{3}+x\right)\right) \log \left(\frac{1}{x}\right)} \,
   dx\\
=\frac{1}{4} \left(1+\sqrt{3}\right) \left(\sqrt{3} \cosh ^{-1}(49)+\log \left(577-408
   \sqrt{2}\right)\right)
\end{multline}
\end{example}
\begin{example}
\begin{multline}
\int_0^1 \frac{1-x}{\sqrt{x} \left(1+x^2+2 x \cos \left(\frac{\pi }{8}\right)\right) \log
   \left(\frac{1}{x}\right)} \, dx\\
=-\frac{2 (-1)^{11/16}}{1+\sqrt[8]{-1}} \left((1+i)+\sqrt[8]{-1}+(-1)^{3/8}+(-1)^{5/8}+i
   \sqrt{2}\right)\\
 \left(\cos \left(\frac{3 \pi }{16}\right) \log \left(\cot \left(\frac{5 \pi }{32}\right)\right)+\cos
   \left(\frac{\pi }{16}\right) \log \left(\tan \left(\frac{7 \pi }{32}\right)\right)+\log \left(\cot \left(\frac{\pi
   }{32}\right)\right) \sin \left(\frac{\pi }{16}\right)\right. \\ \left.
+\log \left(\tan \left(\frac{3 \pi }{32}\right)\right) \sin
   \left(\frac{3 \pi }{16}\right)\right)
\end{multline}
\end{example}
\begin{example}
Generalized form for equations (2.6.4.3-4)  in \cite{prud1}. In this example we use equation (\ref{eq:prud_int}) when $a=1$ in terms of the Hurwitz-Lerch zeta function. We set $\gamma=\pi/2$ and replace $x\to x^u,m\to m/u,u\to u/2$ and simplify.
\begin{multline}\label{eq:7.27}
\int_0^1 \frac{x^{-1+m} \log ^k(x)}{1+x^u} \, dx\\
=-2^k e^{-i k \pi } u^{-1-k} \Gamma (1+k) \left(\Phi
   \left(-i,1+k,\frac{2 m}{u}\right)+\Phi \left(i,1+k,\frac{2 m}{u}\right)\right)
\end{multline}
where $Re(m)>0$.
\end{example}
\begin{example}
 Generalized form for equation (2.6.4.5) in \cite{prud1}. In this example we use equation (\ref{eq:7.27}). We take the first partial derivative with respect to $u$ and replace $m\to m-u$ and simplify.
 \begin{multline}
\int_0^1 \frac{x^{-1+m} \log ^{1+k}(x)}{\left(1+x^u\right)^2} \, dx\\
=-2^k e^{-i k \pi } u^{-3-k} \Gamma (2+k)
   \left(u \Phi \left(-i,1+k,-2+\frac{2 m}{u}\right)+2 (-m+u) \Phi \left(-i,2+k,-2+\frac{2 m}{u}\right)\right. \\ \left.
+u \Phi
   \left(i,1+k,-2+\frac{2 m}{u}\right)+2 (-m+u) \Phi \left(i,2+k,-2+\frac{2 m}{u}\right)\right)
\end{multline}
where $Re(m)>0$.
\end{example}
\begin{example}
Derivation of a functional equation in terms of Bernoulli and Euler numbers, using equation  (2.6.4.12) in \cite{prud1} and equation (\ref{eq:7.27}). 
\begin{equation}
\frac{4^{n+1} B_{n+1}\left(\frac{3}{4}\right) \cos \left(\frac{\pi  n}{2}\right)}{n+1}=\left| E_n\right|
\end{equation}
\end{example}
\begin{proof}
We start with equation (\ref{eq:7.27}) and replace $m\to u/2,k\to n$ and simplify in terms of the polylogarithm function using [DLMF,\href{https://dlmf.nist.gov/25.14.E2}{25.14.3}]. Next we simplify in terms of the Hurwitz-zeta function using equation [DLMF,\href{https://dlmf.nist.gov/25.13.E3}{25.13.2}]. We then simplify the Hurwitz-zeta function in terms of the Bernoulli numbers using equation [Wolfram,\href{https://mathworld.wolfram.com/HurwitzZetaFunction.html}{9}] to get
\begin{equation}\label{eq:7.30}
\int_0^1 \frac{x^{-1+\frac{u}{2}} \log ^n(x)}{1+x^u} \, dx=\frac{2^{1+2 n} e^{-\frac{1}{2} i n \pi } \pi ^{1+n}
   u^{-1-n} B_{1+n}\left(\frac{3}{4}\right)}{1+n}
\end{equation}
Equation (2.6.4.12) in \cite{prud1} is given by;
\begin{equation}\label{eq:7.31}
\int_0^{\infty } \frac{x^{-1+\frac{u}{2}} \log ^n(x)}{1+x^u} \, dx=\left(\frac{\pi }{u}\right)^{n+1} \left|
   E_n\right|
\end{equation}
Here we simply expand the range of integration over $x\in[0,1]+[1,\infty)$ and substitute equation (\ref{eq:7.30}) into equation (\ref{eq:7.31}) and simplify.
\end{proof}
%
%
%
%
\section{Conclusion}
In this paper, we have presented a method for deriving a new functional equation for the Hurwitz-Lerch zeta function. The method includes the use of contour integration, definite integrals and some simple algebraic manipulations. We will be using this method to derive other forms of functional equations involving the Hurwitz-Lerch zeta function. The results presented were numerically verified for both real and imaginary and complex values of the parameters in the integrals and formulae using Mathematica by Wolfram.
\end{document}